\documentclass[amssymb, 11pt]{amsart}
\usepackage{latexsym}

\newlength{\standardunitlength}
\newtheorem{prop}{Proposition}[section]

\newtheorem{lemma}[prop]{Lemma}
\newtheorem{cor}[prop]{Corollary}
\newtheorem{theorem}[prop]{Theorem}

\begin{document}
\newcommand{\ee}{\mathbb{E}}
\newcommand{\pp}{\mathbb{P}}
\newcommand{\aaa}{\mathbb{A}}
\newcommand{\cc}{\mathbb{C}}
\newcommand{\zz}{\mathbb{Z}}

\begin{center}
{\bf Stein's Method and Random Character Ratios}
\end{center}

\begin{center}
{\bf Running head: Stein's Method and Random Character Ratios}
\end{center}

\begin{center}
Version of 8/13/05
\end{center}

\begin{center}
By Jason Fulman
\end{center}

\begin{center}
University of Pittsburgh
\end{center}

\begin{center}
Department of Mathematics, 301 Thackeray Hall, Pittsburgh PA 15260
\end{center}

\begin{center}
email: fulman@math.pitt.edu
\end{center}

{\bf Abstract}: Stein's method is used to prove limit theorems for
random character ratios. Tools are developed for four types of
structures: finite groups, Gelfand pairs, twisted Gelfand pairs, and
association schemes. As one example an error term is obtained for a
central limit theorem of Kerov on the spectrum of the Cayley graph of
the symmetric group generated by i-cycles, or equivalently for the
character ratio of a Plancherel distributed representation on an
i-cycle. Other main examples include an error term for a central limit
theorem of Ivanov on character ratios of random projective
representations of the symmetric group, and a new central limit
theorem for the spectrum of certain graphs whose vertices are the set
of perfect matchings on 2n symbols. The error terms in the resulting
limit theorems are typically $O(n^{-1/4})$ or better. The results
are obtained with remarkably little information: a character formula
for a single representation close to the trivial representation and
estimates on two step transition probabilities of a random
walk. Although the limit theorems stated in this paper are all for the
case of normal approximation, many of the tools developed are quite
general. Indeed, both the construction of an exchangeable pair used
for Stein's method and lemmas computing certain moments are useful for
arbitrary distributional approximation.

\begin{center}
2000 Mathematics Subject Classification: Primary 05E10, Secondary
60C05.
\end{center}

\begin{center}
Key words and phrases: Stein's method, normal approximation, Gelfand
pair, character ratio, symmetric group, Plancherel measure,
association scheme.
\end{center}

\section{Introduction}

	Given a fixed $m \times m$ matrix, it is natural to study
	the distribution of its eigenvalues, where each eigenvalue is
	chosen with probability $\frac{1}{m}$. As a motivating example
	around which discussion can be centered, consider the $n!
	\times n!$ transition matrix for random walk on the symmetric
	group $S_n$, where the generating set consists of all
	$i$-cycles. Diaconis and Shahshahani \cite{DS} proved that the
	eigenvalues of this matrix are the numbers
	$\frac{\chi^{\lambda}_{(i,1^{n-i})}}{dim(\lambda)}$ occurring
	with multiplicity $dim(\lambda)^2$. Here $\lambda$
	parameterizes an irreducible representation of the symmetric
	group, $\chi^{\lambda}_{(i,1^{n-i})}$ is the corresponding
	character value on $i$-cycles, and $dim(\lambda)$ is the
	dimension of the irreducible representation. Since
	$\sum_{|\lambda|=n} dim(\lambda)^2 = n!$, the eigenvalue
	$\frac{\chi^{\lambda}_{(i,1^{n-i})}}{dim(\lambda)}$ is chosen
	with probability $\frac{dim(\lambda)^2}{n!}$.

	The probability measure on irreducible representations of
	$S_n$ which picks the representation corresponding to
	$\lambda$ with probability $\frac{dim(\lambda)^2}{n!}$ is
	known as the Plancherel measure of the symmetric group. Kerov
	\cite{K1} proved that if $i \geq 2$ is fixed, and $\lambda$ is
	random from the Plancherel measure of the symmetric group,
	then the random variable $\frac{\sqrt{{n \choose i} (i-1)!}
	\chi^{\lambda}_{(i,1^{n-i})}}{dim(\lambda)}$ is asymptotically
	normal as $n \rightarrow \infty$. Kerov's proof used the
	method of moments and difficult combinatorics; a beautiful
	exposition of his work is the paper \cite{IO}. Hora \cite{Ho}
	gave another proof of Kerov's result, also using the method of
	moments, but with somewhat simpler combinatorics based on the
	fact that the $k$th moment of the character ratio
	$\frac{\chi^{\lambda}_{(i,1^{n-i})}}{dim(\lambda)}$ is equal
	to the probability that the walk generated by the set of
	i-cycles is at the identity after k steps. In very recent
	work, Sniady \cite{Sn1}, \cite{Sn2} uses the genus
	expansion of random matrix theory to give another method of
	moments proof of Kerov's result.

	A more probabilistic approach to Kerov's result for the case
	$i=2$ was given in \cite{F1}, where Stein's method was used to
	obtain the first error term in Kerov's central limit theorem;
	an error term of $O(n^{-1/4})$ was proved and an error term of
	$O(n^{-1/2})$ was conjectured. This error term was later
	improved to $O(n^{-s})$ for any $s<1/2$ using martingale
	theory \cite{F4}. More recently, a proof of the $O(n^{-1/2})$
	conjecture appears in \cite{SS}, using a new refinement of
	Stein's method. For other proofs of the $O(n^{-1/2})$ bound,
	see \cite{CF} for an argument with a different refinement of
	Stein's method, and \cite{F3} for an argument using
	Bolthausen's variation of Stein's method. All of these
	results, it should be emphasized, were only for the case
	$i=2$. However even in the simple setting of $i=2$, random
	character ratios arise in work on the moduli space of curves
	\cite{EO}.
	
	Before proceeding further, it should be mentioned that
	familiarity with Stein's method is not necessary to read this
	paper. Indeed, Section \ref{stein} gives a very brief
	introduction to normal approximation by Stein's method. It
	presents the bare minimum needed to understand this paper, but
	gives a few pointers to the literature for further reading.

	Section \ref{groups} of this paper generalizes the set-up of
	Kerov's central limit theorem to the case when $G$ is a finite
	group and the generating set consists of a single conjugacy
	class $C$. Quite general results are proved, elucidating our
	early work on this problem \cite{F1} which required different
	information than the current treatment. As a new application,
	it is shown that for any fixed $i \geq 2$, one obtains an
	error term $O(n^{-1/4})$ in Kerov's central limit theorem. The
	approach presented here is appealing because it uses only the
	most elementary ingredients, namely a well known character
	formula for the irreducible representation parameterized by
	$\lambda=(n-1,1)$ on all elements of $S_n$, and estimates on
	the two step transition probabilities of the random walk on
	$S_n$ generated by $C$.

	Section \ref{gelfand} proves theorems in the setting of
	Gelfand pairs. If $G$ is a group and $K$ is a subgroup such
	that the induced representation $1_K^G$ is multiplicity free,
	the pair $(G,K)$ is called a Gelfand pair. Stein's method is
	used to prove a central limit theorem for random spherical
	functions of the pair $(G,K)$. As an application, a new
	central limit theorem with error term $O(n^{-1/4})$ is
	obtained for the spectrum of certain random walks on the set
	of perfect matchings of $2n$ symbols. Equivalently, a central
	limit theorem is obtained for certain statistics under the
	Jack$_{2}$ measure on partitions, which is an interesting
	object \cite{O2}. This application is an exact analog of
	Kerov's central limit theorem, and complements results in the
	paper \cite{F2}, which obtained an analog of Kerov's central
	limit only for $i=2$, but for the entire Jack$_{\alpha}$
	family of measures on partitions, where there is not a known
	corresponding Gelfand pair. As in the group case, only the
	simplest ingredients are needed for the proof: a formula for a
	single spherical function (chosen to be as close to the
	trivial spherical function as possible) and estimates on the
	two step transition probabilities of a random walk.

	Section \ref{twisted} focuses on twisted Gelfand pairs. This
	terminology, introduced in \cite{Stem}, refers to a triple
	$(G,K,\phi)$, where $\phi$ is a linear character of $K$ such
	that $Ind_K^G(\phi)$ is multiplicity free. An error term is
	obtained for a central limit theorem of Ivanov \cite{I} on
	character ratios of random projective representations of the
	symmetric group. There is a close parallel to earlier
	sections, but there are major differences in
	argumentation. One striking example of this is that if one
	attempts a straightforward adaptation of the Markov chain used
	in the construction of an exchangeable pair in Section
	\ref{gelfand}, the resulting ``Markov chain'' can have
	negative transition probabilities. However an interesting
	combinatorial argument shows that only the holding
	``probabilities'' can be negative. Since Stein's method works
	by analyzing how a statistic changes under small
	perturbations, holding probabilities are not so important, and
	by rescaling one obtains a genuine Markov chain and so a
	legitimate construction of an exchangeable pair.

	Section \ref{assoc} develops limit theorems for the spectrum
	of an adjacency matrix of an association scheme. The arguments
	are analogous to those in previous sections. This is not
	surprising since experts in the area will realize that Gelfand
	pairs and association schemes are both generalizations of the
	finite group case. However the perspective and examples are
	quite different, and it would confuse rather than clarify
	matters to present Kerov's central limit theorem as a result
	about Gelfand pairs or association schemes. As an example, the
	Hamming association scheme is treated. One obtains a central
	limit theorem for the spectrum of the Hamming graph, or
	equivalently for values of q-Krawtchouk polynomials.

	Having outlined the contents of this article, it is useful to
give four further reasons why we believe the results to be
interesting. First, as mentioned in the abstract, many of the results
here are applicable for distributional approximations other than
normal approximation. Both the construction of an exchangeable pair
and the computation of certain moments will be quite useful once the
corresponding versions of Stein's method are developed. For an
illustration of this point in the context of the spectrum of
generalizations of the Bernoulli-Laplace diffusion model, see
\cite{CF}. Second, the problem of studying the spectrum of random
walks on $G/K$ where $(G,K)$ is a Gelfand pair, is of ongoing
interest, even in the case when $G$ and $K$ are finite groups. In
particular if $G,K$ are finite classical groups, this leads to
interesting questions in number theory \cite{Te}, and we are working
on using the exchangeable pair constructed in this paper to obtain
information. Third, Plancherel measure (which arises in the group
case), Jack measure (which arises in the Gelfand pair case), and
shifted Plancherel measure (which arise in the twisted Gelfand pair
case), are all objects of interest to researchers in random matrix
theory \cite{AD}, \cite{BO1}, \cite{BOO}, \cite{De}, \cite{J},
\cite{Mat}, \cite{O1}, \cite{O2}, \cite{TW}. As is evident from these
papers, there are many interesting statistics under these measures,
and the method of constructing exchangeable pairs in this paper should
be useful for studying them by Stein's method. Fourth, the examples in
this paper will be a useful testing grounds for results on Stein's
method. For example the refinement of Stein's method in \cite{SS}
arose from trying to obtain an $O(n^{-1/2})$ error term for the $i=2$
case of Kerov's central limit theorem.

	To close the introduction, some deficiencies with this work
	should be noted. First, the method of this paper does not seem
	to be applicable for studying the spectrum of arbitrary
	reversible Markov chains. In all of the examples presented
	here, there is an abundance of symmetry, and the list of
	eigenvalues of the matrix of interest is itself an eigenvector
	for another natural operator. Thus the method of moments,
	which typically doesn't give error terms, may be more
	versatile if one's only purpose is to obtain a limit
	theorem. The examples in the papers \cite{Bi}, \cite{Ho2},
	\cite{HHO} (which are promising candidates for Stein's method
	as well) illustrate the use of the method of moments to study
	the spectrum of random walks. Second, it is not clear that the
	results in this paper are useful for bounding the convergence
	rates of the Markov chains whose spectrum is studied. Such
	bounds, reviewed in Chapter 3 of Diaconis \cite{Di}, usually
	require information about the edge of the spectrum rather than
	the bulk of the spectrum. Third, the error terms, while good,
	are not perfect. As is well known to workers in normal
	approximation, the transition from an $O(n^{-1/4})$ bound to
	an $O(n^{-1/2})$ bound is often quite difficult. We conjecture
	that $O(n^{-1/2})$ bounds hold for Theorems \ref{limgroup},
	\ref{CLTgel}, and \ref{projerror}. To prove this will almost
	certainly require much finer combinatorial information, along
	the lines of that developed in \cite{IO} (which in turn relies
	on \cite{IK}), and even this is only for the special case of
	Plancherel measure of the symmetric group. The point of the
	current paper is to obtain good bounds with minimal
	information and effort.

\section{Stein's Method for Normal Approximation} \label{stein}

	In this section we briefly review Stein's method for normal
	approximation, using the method of exchangeable pairs
	\cite{Stn}. One can also use couplings to prove normal
	approximations by Stein's method (see \cite{Re} for a survey),
	but the exchangeable pairs approach is effective for our
	purposes. For a survey discussing both exchangeable pairs and
	couplings, the paper \cite{RR2} can be consulted.

	Two random variables $W,W'$ on a state space $X$ are called
	exchangeable if for all $w_1,w_2$, $\pp(W=w_1,W'=w_2)$ is
	equal to $\pp(W=w_2,W'=w_1)$. As is typical in probability
	theory, let $\ee(A|B)$ denote the expected value of $A$
	given $B$. The following result of Stein uses an exchangeable
	pair $(W,W')$ to prove a central limit theorem for $W$.

\begin{theorem} \label{steinbound} (\cite{Stn}) Let $(W,W')$ be an
 exchangeable pair of real random variables such that $\ee(W'|W) =
(1-a)W$ with $0< a <1$. Then for all real $x_0$, \begin{eqnarray*} & &
\left| \pp(W \leq x_0) - \frac{1}{\sqrt{2 \pi}} \int_{-\infty}^{x_0}
e^{-\frac{x^2}{2}} dx \right|\\ & \leq &
\frac{\sqrt{Var(\ee[(W'-W)^2|W])}}{a} + (2 \pi)^{- \frac{1}{4}}
\sqrt{\frac{1}{a} \ee|W'-W|^3}. \end{eqnarray*}
\end{theorem}

	In order to apply Theorem \ref{steinbound} to study a
	statistic $W$, one needs an exchangeable pair $(W,W')$. The
	usual way of doing this is to use Markov chain theory. A
	Markov chain $K$ (with chance of going from $x$ to $y$ denoted
	by $K(x,y)$) on a finite set $X$ is called reversible with
	respect to a probability distribution $\pi$ if $\pi(x) K(x,y)
	= \pi(y) K(y,x)$ for all $x,y$. This condition implies that
	$\pi$ is a stationary distribution for $K$. It is
	straightforward to check that if $K$ is reversible with
	respect to $\pi$, then one obtains an exchangeable pair
	$(W,W')$ as follows: choose $x \in X$ from $\pi$, then obtain
	$x'$ by taking one step from $x$ according to $K$, and set
	$(W,W')=(W(x),W(x'))$.

	To apply Theorem \ref{steinbound}, one needs the exchangeable
	pair arising from $K$ to satisfy $\ee(W'|W) = (1-a)W$, and
	also one must obtain useful bounds on $Var(\ee[(W'-W)^2|W])$
	and $\ee|W-W'|^3$. A main contribution of this paper is to
	provide a number of algebraically natural examples where this
	can be done. As the reader will see, the construction of $K$
	in our examples can be subtle.

	A drawback with Theorem \ref{steinbound} is that in many
	problems of interest, it gives a convergence rate of
	$O(n^{-1/4})$ rather than $O(n^{-1/2})$. When $|W'-W|$ is bounded,
	the following variation often gives the correct rate. A
	similar result is in \cite{RR}.

\begin{theorem} \label{rinrot} (\cite{SS}) Let $(W,W')$ be an
 exchangeable pair of real random variables such that $\ee(W'|W) =
(1-a)W$ with $0<a<1$. Suppose that $|W'-W| \leq A$ for some
constant $A$. Then for all real $x_0$, \begin{eqnarray*} & & \left|
\pp(W \leq x_0) - \frac{1}{\sqrt{2 \pi}} \int_{-\infty}^{x_0}
e^{-\frac{x^2}{2}} dx \right|\\ & \leq &
\frac{\sqrt{Var(\ee[(W'-W)^2|W])}}{a} + .41 \frac{A^3}{a} + 1.5
A. \end{eqnarray*}
\end{theorem} In many of our examples, $|W'-W|$ is unbounded. Thus as
	mentioned in the introduction, significantly more work may be
	needed to go beyond the $O(n^{-1/4})$ rate.

	The following general lemmas are helpful in working with the
	error bounds in either Theorem \ref{steinbound} or Theorem
	\ref{rinrot}.

\begin{lemma} \label{var} Let $(W,W')$ be an exchangeable pair of random variables such that $\ee(W'|W) = (1-a) W$ and $\ee(W^2)=1$. Then $\ee(W'-W)^2 = 2a$. 
\end{lemma}

\begin{proof} Since $W$ and $W'$ have the same distribution, \begin{eqnarray*}
 \ee(W'-W)^2 & = & \ee(\ee(W'-W)^2|W)\\ & = & \ee((W')^2) + \ee(W^2) -
2 \ee(W \ee(W'|W))\\ & = & 2 \ee(W^2) - 2 \ee(W \ee(W'|W))\\ & = & 2
 \ee(W^2) - 2(1-a) \ee(W^2)\\ & = & 2a. \end{eqnarray*} \end{proof}

	Lemma \ref{majorize} is a well known inequality (already used
	in the monograph \cite{Stn}) and useful because often the
	right hand side is easier to compute or bound than the left
	hand side. To make this paper as self-contained as possible,
	we include a proof. Here $x$ is an element of the state space
	$X$.

\begin{lemma} \label{majorize} \[ Var(\ee[(W'-W)^2|W]) \leq Var (\ee[(W'-W)^2|x]).\] \end{lemma}	

\begin{proof} Jensen's
 inequality states that if $g$ is a convex function, and $Z$ a random
 variable, then $g(\ee(Z)) \leq \ee(g(Z))$. There is also a
 conditional version of Jensen's inequality (Section 4.1 of \cite{Du})
 which states that for any $\sigma$ subalgebra ${\it F}$ of the
 $\sigma$-algebra of all subsets of $X$, \[ \ee(g(\ee(Z|{\it F})))
 \leq \ee(g(Z)).\] The lemma follows by setting $g(t)=t^2$,
 $Z=\ee((W'-W)^2|x)$, and letting ${\it F}$ be the $\sigma$-algebra
 generated by the level sets of $W$. \end{proof}

\section{Finite groups} \label{groups}

	This section uses Stein's method to study the spectrum of
	random walk on a finite group $G$, where the generating set is
	a conjugacy class $C$ which satisfies $C=C^{-1}$. As will be
	explained in Subsection \ref{rtgroup}, by \cite{DS} this is
	equivalent to studying the distribution of the character
	ratio $\frac{\chi^{\lambda}(C)}{dim(\lambda)}$ where
	$\lambda$ is chosen from the Plancherel measure of the group
	$G$.

	The organization of this section is as follows. Subsection
	\ref{rtgroup} recalls the necessary background from
	representation theory. Subsection \ref{groupCLT} then defines
	a Markov chain on the set of irreducible representations of
	$G$, and uses it to construct an exchangeable pair. It then
	shows that the moment computations needed for Stein's method
	can be carried out provided that one knows the character
	values of a single nontrivial irreducible representation of
	$G$ (which to be useful should be close to the trivial
	representation) and that one can estimate the two-step
	transition probabilities of the random walk on $G$ generated
	by $C$. This leads to a general central limit
	theorem. Subsection \ref{Sn} applies the theory to the
	symmetric group $S_n$ with $C$ the conjugacy class of
	$i$-cycles, where $i$ is fixed and $n$ is large.

\subsection{Background from representation theory} \label{rtgroup}

	We recall facts from the representation theory of finite
	groups, referring the reader to \cite{Sa} or \cite{Se} for
	more details. In what follows, $\chi$ denotes a character of
	the finite group $G$, $dim(\chi)$ is the dimension of the
	corresponding representation, and $Irr(G)$ is the set of all
	irreducible characters of $G$. Also $\overline{z}$ denotes
	complex conjugate of a number $z$.

\begin{lemma} \label{inverse} Let $\chi$ be an irreducible representation of $G$. Then $\chi(C^{-1}) = \overline{\chi(C)}$. Thus if $C=C^{-1}$, then $\chi(C)$ is real. \end{lemma}

	Next we recall the orthogonality relations for irreducible
	characters of $G$.

\begin{lemma} \label{orth0} Let $\nu$ and $\chi$ be irreducible characters of a finite group $G$. Then  \[ \frac{1}{|G|} \sum_{g \in G} \nu(g)
\overline{\chi(g)} = \delta_{\nu,\chi}. \] \end{lemma}

\begin{lemma} \label{orth1} Let $C$ be the conjugacy class of $G$ containing
the element $g$. Then for $g \in G$, \[ \sum_{\chi \in Irr(G)} \chi(g)
\overline {\chi(h)} \] is equal to $\frac{|G|}{|C|}$ if $h,g$ are
conjugate and is $0$ otherwise.
\end{lemma}

	Lemma \ref{countsol} while known, is perhaps not well known,
	and since analogous results will be needed in later sections,
	a proof along the lines of one in \cite{HSS} is included.

\begin{lemma} \label{countsol} Let $G$ be a finite group with conjugacy classes $C_1, \cdots, C_t$. Let $C_k$ be the
 conjugacy class of an element $w \in G$. Then the number of
 $m$-tuples $(g_1,\cdots,g_m) \in G^m$ such that $g_j \in C_{i_j}$ and
 $g_1 \cdots g_m=w$ is \[ \prod_{j=1}^m |C_{i_j}|
 \sum_{\chi \in Irr(G)} \frac{dim(\chi)^2}{|G|}
 \frac{\chi(C_{i_1})}{dim(\chi)} \cdots \frac{
 \chi(C_{i_m})}{dim(\chi)} \overline{\frac{\chi(C_k)}{dim(\chi)}} .\]
 \end{lemma}

\begin{proof} Identify each class $C_i$ with its
 corresponding class sum in the complex group algebra $\cc G$. If
 $\chi^{(1)},\cdots,\chi^{(t)}$ are the irreducible complex characters
 of $G$, then the elements \[ E_s =
 \frac{dim(\chi^{(s)})}{|G|} \sum_{j=1}^t \overline{\chi_j}^{(s)} C_j
 \ \ \ (1 \leq s \leq t) \] are a complete set of orthogonal
 idempotents for the center of $\cc G$, where $\chi_j^{(s)}$ denotes
 the value of $\chi^{(s)}$ at any $g \in C_j$. Lemma
 \ref{orth1} implies that \[ C_j = |C_j| \sum_{s=1}^t
 \frac{\chi_j^{(s)}}{dim(\chi^{(s)})} E_s.\] Since the $E_j$'s are
 orthogonal idempotents (i.e. $E_r E_s = \delta_{r,s} E_r$), it
 follows that \begin{eqnarray*} C_{i_1} C_{i_2} \cdots C_{i_m} & = & |C_{i_1}| \cdots
 |C_{i_m}| \sum_{s=1}^t \frac{\chi_{i_1}^{(s)} \cdots \chi_{i_m}^{(s)}}{dim(\chi^{(s)})^m}
 E_s\\ & = & \frac{|C_{i_1}| \cdots |C_{i_m}|}{|G|} \sum_{k=1}^t C_k
 \sum_{s=1}^t \frac{ \chi_{i_1}^{(s)} \cdots \chi_{i_m}^{(s)}  \overline{\chi_k}^{(s)} }{dim(\chi^{s})^{m-1}}
, \end{eqnarray*} as desired. \end{proof}

	As a corollary one obtains the following result.

\begin{cor}
 \label{specgroup} (\cite{DS}) Suppose that $C$ is a conjugacy class
 satisfying $C=C^{-1}$. Then the eigenvalues of the random walk on $G$
 with generating set $C$ are indexed by $\chi \in Irr(G)$ and are the
 numbers $\frac{\chi(C)}{dim(\chi)}$, occurring with multiplicity
 $dim(\chi)^2$. \end{cor}

\begin{proof} If $M$ is the $|G| \times |G|$ transition matrix for the random walk, the chance of being at the identity after $k$ steps is the trace of $M^k$ divided by $|G|$. Thus Lemma \ref{countsol} implies that for all $k \geq 0$, the trace of $M^k$ is equal to \[ \sum_{\chi \in Irr(G)} dim(\chi)^2
 \left( \frac{\chi(C)}{dim(\chi)} \right)^k .\] Since this holds for
 all $k \geq 0$, the result follows. \end{proof}

	As mentioned in the introduction, the Plancherel measure of
	$G$ is the probability measure on $Irr(G)$ which chooses each
	$\chi$ with probability $\frac{dim(\chi)^2}{|G|}$. So
	Corollary \ref{specgroup} says that the eigenvalues of the
	random walk on $G$ generated by $C$ are the ``character
	ratios'' $\frac{\chi(C)}{dim(\chi)}$ occurring with
	multiplicity proportional to the Plancherel probability of
	$\chi$.

\subsection{Central limit theorems for character ratios} \label{groupCLT}
	
	The goal of this subsection is to prove a central limit
	theorem for the random variable $W$ defined by $W(\lambda) =
	\frac{|C|^{1/2} \chi^{\lambda}(C)}{dim(\lambda)}$, where $C$
	is a fixed conjugacy class such that $C=C^{-1}$ and $\lambda$
	is random from the Plancherel measure of $G$.

	For this purpose, it is useful to construct a Markov chain on
	the set of irreducible representations of $G$ as
	follows. First, fix a non-trivial irreducible representation
	$\tau$ whose character is real valued. This gives a Markov chain
	$L_{\tau}$ by defining the probability of transitioning from
	$\lambda$ to $\rho$ as \[ L_{\tau}(\rho,\lambda) :=
	\frac{dim(\rho)}{dim(\lambda) dim(\tau)} \frac{1}{|G|}
	\sum_{g} \chi^{\lambda}(g) \chi^{\tau}(g)
	\overline{\chi^{\rho}(g)}. \] It is worth remarking that since
	Stein's method is in the spirit of Taylor approximation,
	$\tau$ should typically be chosen to be as close as possible
	to the trivial representation.

	Lemma \ref{gchain} verifies that $L_{\tau}$ is a Markov chain
which is reversible with respect to Plancherel measure.

\begin{lemma} \label{gchain} The transition probabilities of $L_{\tau}$ are real and non-negative and sum to 1. Moreover the chain $L_{\tau}$ is reversible with respect to the Plancherel measure of $G$.
\end{lemma}

\begin{proof} The transition probabilities of $L_{\tau}$ are real and non-negative since \[  \frac{1}{|G|} \sum_{g} \chi^{\lambda}(g) \chi^{\tau}(g)
	\overline{\chi^{\rho}(g)} \] is the multiplicity of $\rho$ in
	the tensor product of $\lambda$ and $\tau$. Letting $id$
	denote the identity, it follows from Lemma \ref{orth1} that \[
	\sum_{\rho} L_{\tau}(\lambda,\rho) = \frac{1}{dim(\lambda)
	dim(\tau) |G|} \sum_g \chi^{\lambda}(g) \chi^{\tau}(g)
	\sum_{\rho} \chi^{\rho}(id) \overline{\chi^{\rho}(g)} = 1.\]
	For the reversibility assertion, the fact that $\chi^{\tau}$
	and the transition probabilities of $L_{\tau}$ are both real
	valued implies that \[ \frac{dim(\lambda)^2}{|G|}
	L_{\tau}(\lambda,\rho) = \frac{dim(\rho)^2}{|G|}
	L_{\tau}(\rho,\lambda). \] \end{proof}

	An exchangeable pair $(W,W')$ is now constructed from the
	chain $L_{\tau}$ in the standard way. First choose $\lambda$
	from the Plancherel measure of $G$, then choose $\rho$ with
	probability $L_{\tau}(\lambda,\rho)$, and finally let $(W,W')
	= (W(\lambda),W(\rho))$. The remaining results in this
	subsection show that the exchangeable pair $(W,W')$ has
	desirable properties.

\begin{lemma} \label{steinsat} $\ee(W'|W) = \left(  \frac{\chi^{\tau}(C)}{dim(\tau)} \right) W$. \end{lemma} 

\begin{proof} From the definition of $W'$,
\begin{eqnarray*}
\ee(W'|\lambda) & = & |C|^{1/2} \sum_{\rho} \frac{dim(\rho)}{dim(\lambda) dim(\tau)} \frac{1}{|G|} \sum_g \chi^{\lambda}(g) \chi^{\tau}(g) \overline{\chi^{\rho}(g)} \frac{\chi^{\rho}(C)}{dim(\rho)}\\
& = & \frac{|C|^{1/2}}{|G|} \sum_{g} \frac{\chi^{\lambda}(g)}{dim(\lambda)} \frac{\chi^{\tau}(g)}{dim(\tau)} \sum_{\rho} \chi^{\rho}(C) \overline{\chi^{\rho}(g)}\\
& = &  \left(  \frac{\chi^{\tau}(C)}{dim(\tau)} \right) W(\lambda), \end{eqnarray*} where the last step is Lemma \ref{orth1}. The result follows since this depends on $\lambda$ only through $W$. \end{proof}

	Corollary \ref{record} is not needed in what follows, but is
	worth recording.

\begin{cor}
 \label{record} The eigenvalues of $L_{\tau}$ are
 $\frac{\chi^{\tau}(C)}{dim(\tau)}$ as $C$ ranges over conjugacy
 classes of $G$. The functions $\psi_C(\lambda) = \frac{|C|^{1/2}
 \chi^{\lambda}(C)}{dim(\lambda)}$ are a basis of eigenvectors of
 $L_{\tau}$, orthonormal with respect to the inner product \[ \langle
 f_1,f_2 \rangle = \sum_{\lambda} f_1(\lambda) \overline{f_2(\lambda)}
 \frac{dim(\lambda)^2}{|G|}.\]
\end{cor}

\begin{proof} The proof of Lemma \ref{steinsat} shows that $\psi_C$ is an eigenvector of $L_{\tau}$ with eigenvalue $\frac{\chi^{\tau}(C)}{dim(\tau)}$. The orthonormality assertion follows from Lemma \ref{orth1}, and the basis assertion follows since the number of conjugacy classes of $G$ is equal to the number of irreducible representations of $G$. \end{proof} 

\begin{lemma} \label{immcor} $\ee(W'-W)^2 = 2 \left(1 - \frac{\chi^{\tau}(C)}{dim(\tau)} \right)$. \end{lemma}

\begin{proof} This is immediate from Lemma \ref{var} and \ref{steinsat}. \end{proof}

	For the remainder of this subsection, if $K$ is a conjugacy
	class of $G$, $p_m(K)$ will denote the probability that the
	random walk generated by $C$ is in $K$ after $m$ steps.

\begin{lemma} \label{pre2} \[ \ee((W')^2|\lambda) = \frac{|C|}{dim(\lambda) dim(\tau)} \sum_K p_2(K) \chi^{\lambda}(K) \chi^{\tau}(K) \] where the sum is over all conjugacy classes $K$ of $G$.     
\end{lemma}

\begin{proof} \begin{eqnarray*} & & \ee((W')^2|\lambda)\\ & = & \frac{1}{ dim(\lambda) dim(\tau) |G|} \sum_{\rho} dim(\rho) \sum_g \chi^{\lambda}(g) \chi^{\tau}(g) \overline{\chi^{\rho}(g)} \left(\frac{|C|^{1/2} \chi^{\rho}(C)}{dim(\rho)} \right)^2 \\
& = & \frac{|C|}{ dim(\lambda) dim(\tau) |G|} \sum_g \chi^{\lambda}(g)
\chi^{\tau}(g) \sum_{\rho} dim(\rho) \overline{\chi^{\rho}(g)} \left(
\frac{\chi^{\rho}(C)}{dim(\rho)} \right)^2. \end{eqnarray*} The result
now follows from Lemma \ref{countsol}. \end{proof}

	Lemma \ref{big1} writes $Var([\ee(W'-W)^2|\lambda])$ as a sum
	of positive quantities.

\begin{lemma} \label{big1} \[
 Var([\ee(W'-W)^2|\lambda]) = |C|^2 \sum_{ K \neq id}
 \frac{p_2(K)^2}{|K|} \left( \frac{\chi^{\tau}(K)}{dim(\tau)} + 1 -
 \frac{2 \chi^{\tau}(C)}{dim(\tau)} \right)^2  , \] where $K$
 ranges over all non-identity conjugacy classes of $G$. \end{lemma}

\begin{proof} By Lemmas \ref{pre2}, \ref{immcor}, and \ref{steinsat}, \begin{eqnarray*} & & Var([\ee(W'-W)^2|\lambda])\\ & = & \ee \left( \ee((W'-W)^2|\lambda)^2 \right) - 4 \left( 1 - \frac{\chi^{\tau}(C)}{dim(\tau)} \right)^2\\
& = & |C|^2 \ee \left[ \frac{ \sum_K p_2(K)
\chi^{\lambda}(K) \chi^{\tau}(K)   }{dim(\lambda) dim(\tau)}
 + \left(1 - \frac{2 \chi^{\tau}(C)}{dim(\tau)} \right) \left(\frac{\chi^{\lambda}(C)}{dim(\lambda)} \right)^2 \right]^2\\
& &  - 4  \left( 1 - \frac{\chi^{\tau}(C)}{dim(\tau)} \right)^2 \\ & = & T_1 + T_2 + T_3 - 4  \left( 1 - \frac{\chi^{\tau}(C)}{dim(\tau)} \right)^2 \end{eqnarray*} where
\[ T_1 = \frac{|C|^2}{ dim(\tau)^2} \ee \left( \sum_K \frac{p_2(K) \chi^{\lambda}(K) \chi^{\tau}(K)}{dim(\lambda)} \right)^2 \]

\[ T_2 = \frac{2 |C|^2}{ dim(\tau)}
 \left(1-\frac{2 \chi^{\tau}(C)}{dim(\tau)} \right) \sum_K p_2(K)
 \chi^{\tau}(K) \ee \left[ \frac{\chi^{\lambda}(K)}{dim(\lambda)}
\left( \frac{\chi^{\lambda}(C)}{dim(\lambda)} \right)^2 \right] .\]

\[ T_3 =  |C|^2 \left( 1-\frac{2 \chi^{\tau}(C)}{dim(\tau)} \right)^2 \ee \left( \frac{\chi^{\lambda}(C)}{dim(\lambda)} \right)^4 .\]

	Fortunately, these terms can be simplified. By Lemma
	\ref{pre2}, the expression inside parentheses in $T_1$ is
	real. Since $\chi^{\tau}$ is real valued, Lemma \ref{orth1}
	implies that \[ T_1 = |C|^2 \sum_K \frac{p_2(K)^2}{|K|}
	\frac{\chi^{\tau}(K)^2}{dim(\tau)^2}.\] Since
	$\chi^{\lambda}(C)$ is real, Lemma \ref{countsol} implies that
	\[ T_2 = 2 |C|^2 \left( 1-\frac{2 \chi^{\tau}(C)}{dim(\tau)}
	\right) \sum_K \frac{p_2(K)^2}{|K|} \frac{\chi^{\tau}(K)}
	{dim(\tau)} \] and since $C=C^{-1}$, \begin{eqnarray*} T_3 & =
	&  |C|^2 \left( 1-\frac{2 \chi^{\tau}(C)}{dim(\tau)} \right)^2 p_4(id)\\ & =
	&  |C|^2 \left( 1-\frac{2 \chi^{\tau}(C)}{dim(\tau)} \right)^2 \sum_K
	\frac{p_2(K)^2}{|K|}. \end{eqnarray*} Thus one can write \[
	T_1+T_2+T_3 = |C|^2 \sum_K \frac{p_2(K)^2}{|K|} \left(
	\frac{\chi^{\tau}(K)}{dim(\tau)} + 1 - 
	\frac{2 \chi^{\tau}(C)}{dim(\tau)} \right)^2 \] and the result
	follows since $p_2(id) = \frac{1}{|K|}$. \end{proof}

\begin{lemma} \label{mom1} Let $k$ be a positive integer.
\begin{enumerate}
\item $\ee(W'-W)^k= |C|^{k/2} \sum_{m=0}^k  (-1)^{k-m} {k \choose m} \sum_K \frac{\chi^{\tau}(K)}{dim(\tau)} \frac{p_m(K) p_{k-m}(K)}{|K|}$.
\item $\ee(W'-W)^4$ is equal to \[ |C|^2 
\sum_K \left[  8 \left(1 - \frac{\chi^{\tau}(C)}{dim(\tau)} \right) - 6 \left( 1-\frac{\chi^{\tau}(K)}{dim(\tau)} \right) \right] \frac{p_2(K)^2}{|K|}.\]
\end{enumerate}
\end{lemma}  

\begin{proof} For the first assertion, note that
\begin{eqnarray*}  \ee((W'-W)^k|\lambda)
 & = & \frac{|C|^{k/2}}{dim(\lambda) dim(\tau)}
\sum_{\rho} \frac{dim(\rho)}{|G|}
\sum_g \chi^{\lambda}(g) \chi^{\tau}(g) \overline{\chi^{\rho}(g)}\\ & & \cdot \sum_{m=0}^k
(-1)^{k-m} {k \choose m} \left(\frac{\chi^{\rho}(C)}{dim(\rho)}
\right)^m \left( \frac{\chi^{\lambda}(C)}{dim(\lambda)} \right)^{k-m}
\\ & = & \frac{|C|^{k/2}}{dim(\lambda) dim(\tau)} \sum_{m=0}^k
(-1)^{k-m} {k \choose m} \left( \frac{\chi^{\lambda}(C)}{dim(\lambda)}
\right)^{k-m}\\ & & \cdot \sum_g \chi^{\tau}(g) \chi^{\lambda}(g)
\sum_{\rho} \frac{dim(\rho)}{|G|} \left( \frac{\chi^{\rho}(C)}{dim(\rho)}
\right)^m \overline{\chi^{\rho}(g)} \\
& = & \frac{|C|^{k/2}}{dim(\lambda) dim(\tau)}
\sum_{m=0}^k (-1)^{k-m} {k \choose m} \left(
\frac{\chi^{\lambda}(C)}{dim(\lambda)} \right)^{k-m}\\ & & \cdot \sum_K \chi^{\
\tau}(K) \chi^{\lambda}(K) p_m(K), \end{eqnarray*} where the final equality is by Lemma \ref{countsol}. Thus $ \ee((W'-W)^k)$ is equal to   \begin{eqnarray*}
 &  & \ee(\ee((W'-W)^k|\lambda))\\
& = &  \frac{|C|^{k/2}}{dim(\tau)}
\sum_{m=0}^k (-1)^{k-m} {k \choose m} \\ & & 
\cdot \sum_K p_m(K) \chi^{\tau}(K) \sum_{\lambda} \frac{dim(\lambda)^2}{|G|} \frac{\chi^{\lambda}(K)}{dim(\lambda)} \left( \frac{\chi^{\lambda}(C)}{dim(\lambda)} \right)^{k-m}. \end{eqnarray*} The first assertion now follows from Lemma \ref{countsol} and the fact that $\chi^{\lambda}(C)$ is real for all $\lambda$.

	For the second assertion, note by the first assertion that \[
	\ee(W'-W)^4 = |C|^2 \sum_{m=0}^4 (-1)^{m} {4 \choose m} \sum_K
	\frac{\chi^{\tau}(K)}{dim(\tau)} \frac{p_m(K) p_{4-m}(K)}{|K|}
	.\] If $\tau$ is the trivial representation, then $W'=W$,
	which implies that \[ 0 = |C|^2 \sum_{m=0}^4 (-1)^{m} {4
	\choose m} \sum_K \frac{p_m(K) p_{4-m}(K)}{|K|}. \] Thus for
	general $\tau$, \[ \ee(W'-W)^4 = - |C|^2 \sum_{m=0}^4 (-1)^{m}
	{4 \choose m} \sum_K \left( 1 -
	\frac{\chi^{\tau}(K)}{dim(\tau)} \right) \frac{p_m(K)
	p_{4-m}(K)}{|K|}.\] Observe that the $m=0,4$ terms in this sum
	vanish, since the only contribution could come from the
	identity, which contributes 0. The $m=2$ term is \[ -6 |C|^2
	\sum_K \left( 1-\frac{\chi^{\tau}(K)}{dim(\tau)} \right)
	\frac{p_2(K)^2}{|K|} .\] The $m=1,3$ terms are equal and
	together contribute \begin{eqnarray*} & & 8 |C|^2 \sum_K \left(
	1-\frac{\chi^{\tau}(K)}{dim(\tau)} \right) \frac{p_1(K)
	p_3(K)}{|K|}\\ & = & 8 |C| \left(1 -
	\frac{\chi^{\tau}(C)}{dim(\tau)} \right) p_3(C)\\ & = & 8
	|C|^2 \left(1 - \frac{\chi^{\tau}(C)}{dim(\tau)} \right)
	p_4(id)\\ & = & 8 |C|^2 \left(1 -
	\frac{\chi^{\tau}(C)}{dim(\tau)} \right) \sum_K
	\frac{p_2(K)^2}{|K|}, \end{eqnarray*} where $id$ is the
	identity and the last two equalities can be seen either
	directly or from Lemma \ref{countsol}. This completes the
	proof of the second assertion. \end{proof}

	Putting the pieces together, one obtains the following theorem.

\begin{theorem} \label{main1} Let
 $C$ be a conjugacy class of a finite group $G$ such that $C=C^{-1}$
 and fix a nontrivial irreducible representation $\tau$ of $G$ whose character
 is real valued. Let $\lambda$ be a random irreducible representation,
 chosen from the Plancherel measure of $G$. Let $W=\frac{|C|^{1/2}
 \chi^{\lambda}(C)}{dim(\lambda)}$. Then for all real $x_0$,
 \begin{eqnarray*} && \left| \pp(W \leq x_0) - \frac{1}{\sqrt{2 \pi}}
 \int_{-\infty}^{x_0} e^{-\frac{x^2}{2}} dx \right|\\ & \leq &
 \frac{|C|}{a} \sqrt{\sum_{ K \neq id} \frac{p_2(K)^2}{|K|} \left(
 \frac{\chi^{\tau}(K)}{dim(\tau)} + 2a -1 \right)^2}\\ & & +
  \left[ \frac{|C|^2}{\pi} \sum_K \left(8 - \frac{6}{a} \left(
 1-\frac{\chi^{\tau}(K)}{dim(\tau)} \right) \right)
 \frac{p_2(K)^2}{|K|} \right]^{1/4} , \end{eqnarray*} where
 $a=1-\frac{\chi^{\tau}(C)}{dim(\tau)}$. \end{theorem}

\begin{proof}
 One applies Theorem \ref{steinbound} to the exchangeable pair
 $(W,W')$ of this subsection. By Lemmas \ref{majorize} and \ref{big1},
 the first term in Theorem \ref{steinbound} gives the first term in
 the theorem. To upper bound the second term in Theorem
 \ref{steinbound}, note by the Cauchy-Schwarz inequality that \[
 \ee|W'-W|^3 \leq \sqrt{ \ee(W'-W)^2 \ee(W'-W)^4}.\] Now use Lemma
 \ref{immcor} and part 2 of Lemma \ref{mom1}. \end{proof}

\subsection{Example: Cayley Graphs of the Symmetric Group} \label{Sn}

	This subsection applies the theory of Subsection
	\ref{groupCLT} to the case where the group is $S_n$ and $C$ is
	the conjugacy class of $i$-cycles.

	It is useful to recall some facts about the symmetric
	group. Since $K=K^{-1}$ for all conjugacy classes $K$ of
	$S_n$, Lemma \ref{inverse} implies that all irreducible
	characters of $S_n$ are real valued. Also it is elementary
	that $|K| = \frac{n!}{\prod_j j^{m_j} m_j!}$, where $m_j$ is
	the number of cycles of length $j$ of an element of $K$.

	In order to upper bound the error terms in Theorem \ref{main1}, the
	following estimate on the two step transition probabilities of
	the random walk generated by $C$ will be useful. 

\begin{lemma} \label{term1} Let $C$ be the conjugacy class of cycles of length i of the symmetric group $S_n$. Then for $i$ fixed and $n \geq 2i$, $\frac{p_2(K)^2}{|K|}$ is equal to 
\begin{enumerate}
\item $i^2 n^{-2i} + O(n^{-2i-1})$ if $K$ is the identity conjugacy
class.
\item $2 i^2 n^{-2i} + O(n^{-2i-1})$ if $K$ is the conjugacy class consisting of two cycles of length i.
\item $O(n^{-2i-1})$ otherwise.
\end{enumerate}
\end{lemma}

\begin{proof} The first assertion is clear since if $K$ is the identity
 class, $p_2(K)= \frac{1}{|C|}$. For the second assertion, note that
 $p_2(K)=1+O(n^{-1})$ since the only way that a product of two i-cycles is
 not in $K$ is if the i-cycles have a symbol in common. 

	For the third assertion, there are two cases. The first case
	is that $K$ has exactly $n-2i$ fixed points. Then $|K|$ is at
	least $c_i n^{2i}$ where $c_i$ is a constant depending on
	$i$. Since $K$ is not the class consisting of exactly two
	i-cycles, $p_2(K)=O(n^{-1})$ by the second assertion, proving
	the third assertion in this case. The second case is that $K$ has
	$n-2i+r$ fixed points, where $1 \leq r < 2i$. Then $|K|$ is
	at least $c_{i} n^{2i-r}$ where $c_i$ is a constant depending
	on $i$. However $p_2(K) = O(n^{-r})$, since in the product of
	the two i-cycles, there are r symbols moved by the first
	i-cycle each of which is mapped back to itself by the second
	i-cycle. Thus in this case also, the third assertion is
	proved. \end{proof}

\begin{theorem} \label{limgroup} Let $C$ be the conjugacy
 class of $i$-cycles in $S_n$, with $i \geq 2$. Choosing $\lambda$
 from the Plancherel measure of the symmetric group, define a random
 variable $W=\frac{\sqrt{{n \choose i} (i-1)!} \chi^{\lambda}(C)}{dim(\lambda)}$. Then
 there is a constant $A_i$ such that for all real $x_0$, \[ \left|
 \pp(W \leq x_0) - \frac{1}{\sqrt{2 \pi}} \int_{-\infty}^{x_0}
 e^{-\frac{x^2}{2}} dx \right| \leq A_i n^{-1/4} .\] \end{theorem}

\begin{proof} One applies Theorem \ref{main1}, choosing $\tau$ to be the irreducible representation corresponding to $\tau=(n-1,1)$, whose
 character value is the number of fixed points $- 1$. Then for any
 non-identity conjugacy class $K$, \[ \frac{p_2(K)^2}{|K|} \left(
 \frac{\chi^{\tau}(K)}{dim(\tau)} + 2a - 1 \right)^2 = O(n^{-2i-3}).\] This follows from Lemma \ref{term1} and the fact that $ \frac{\chi^{\tau}(K)}{dim(\tau)}+ 2a -1$ vanishes if $K$ has $n-2i$ fixed points
 and is $O(n^{-1})$ for any other class $K$ with $p_2(K) \neq 0$,
 since such a class has at least $n-2i$ fixed points.  Since $|C|=
 \frac{n!}{(n-i)! i}$ and $a=\frac{i}{n-1}$, it follows that the first error term in Theorem \ref{main1} is at
 most $A_i' n^{-1/2}$, where $A_i'$ is a constant depending on $i$.

To bound the second error term in Theorem \ref{main1}, observe that
 Lemma \ref{term1} implies that \[ \sum_K \left(8 - \frac{6}{a} \left(
 1-\frac{\chi^{\tau}(K)}{dim(\tau)} \right) \right)
 \frac{p_2(K)^2}{|K|} \] is $O(n^{-2i-1})$, since the $n^{-2i}$
 contributions from the identity class and the class of two i-cycles
 cancel. Hence the second error term is at most $A_i'' n^{-1/4}$ where
 $A_i''$ is a constant depending only on $i$. \end{proof}

\section{Gelfand pairs} \label{gelfand}

	If $G$ is a finite group and $K$ is a subgroup of $G$ such
	that the induced representation $1_K^G$ is multiplicity free,
	then the pair $(G,K)$ is called a Gelfand pair. This section
	shows that the results of Section \ref{groups} have an
	analog for Gelfand pairs.

	Subsection \ref{rtgelfand} discusses the representation theory
	of Gelfand pairs. Subsection \ref{gelfandCLT} derives a
	general central limit theorem for random spherical functions
	of a Gelfand pair. Subsection \ref{hypercube} illustrates the
	theory on a toy example, proving a central limit theorem for
	the spectrum of the hypercube. A more serious example is
	considered in Subsection \ref{matching}, which obtains a
	new central limit theorem for the spectrum of certain graphs whose
	vertices are the perfect matchings of 2n symbols.

\subsection{Background from representation theory} \label{rtgelfand}

        We discuss some facts about the representation theory of
        Gelfand pairs. A useful reference is Chapter 7 of \cite{Mac}.

	The induced representation $1_K^G$ decomposes in a
        multiplicity free way as $\oplus_{r=0}^s V_r$, where $V_0$
        denotes the trivial representation of $G$. Let $d_r$ denote
        the dimension of $V_r$. For $0 \leq r \leq s$, let $\omega_r$
        denote the corresponding spherical function on $G$, defined by \[
        \omega_r(x) = \frac{1}{|K|} \sum_{k \in K} \chi^{(r)} (x^{-1}
        k) \] where $\chi^{(r)}$ is the character of $V_r$. The
        functions $\omega_r$ are a basis of the space of functions on
        $G$ which are constant on the double cosets of $K$ in $G$. Let
        $K_0,\cdots,K_s$ denote the double cosets of $K$ in $G$ and
        let $g_0,\cdots,g_s$ be corresponding double coset
        representatives, so that $K_i = K g_i K$. It is convenient to
        take $g_0$ to be the identity element of $G$. From page 389 of
        \cite{Mac}, one has that $\omega_i(g_0)=1$ for all $i$.

	It is useful to recall facts which are analogs of those in
	Subsection \ref{rtgroup}.

\begin{lemma} \label{invgelfand} (\cite{Mac}, page 389) Let $\omega$ be a spherical function of the Gelfand pair $(G,K)$. Then $\omega(x^{-1})= \overline{\omega(x)}$.
\end{lemma}

	The following two orthogonality relations are also useful.

\begin{lemma} \label{orthog1} (\cite{Mac}, page 389) For $0 \leq i,j \leq s$, \[
\frac{d_i}{|G|} \sum_{r=0}^s |K_r| \omega_i(g_r)
\overline{\omega_j(g_r)} = \delta_{i,j}.\] \end{lemma}

\begin{lemma} \label{orthog2} For $0 \leq r,t \leq s$, \[ \sum_{i=0}^s d_i
 \omega_i(g_r) \overline{\omega_i(g_t)} = \delta_{r,t}
\frac{|G|}{|K_r|}.\]
\end{lemma}

\begin{proof} Consider the $s \times s$ matrix whose entry in the ith column and rth row is $\sqrt{ \frac{d_i |K_r|}{|G|}} \omega_i(g_r)$. By Lemma \ref{orthog1}, the columns of this matrix are orthonormal. Hence so are its rows, proving the lemma. \end{proof} 

	If $P$ is a $K$ biinvariant probability on $G$ (i.e. constant
	on double cosets of $K$ in $G$), let $p_m(K_r)$ denote the
	probability that the m-fold convolution of $P$ assigns to the
	double coset $K g_r K$. Lemma \ref{fourier} is an analog of
	Lemma \ref{countsol} and could be proved along similar lines,
	as in \cite{HSS}. Instead, we use the language of Fourier
	analysis, as developed on page 395 of \cite{Mac}.

\begin{lemma} \label{fourier} Let
 $P$ be the $K$ biinvariant probability on $G$ which associates mass
 $1$ to the double coset $K_u$ and mass 0 to all other $K$ double
 cosets of $G$. Then for $0 \leq r \leq s$, \[ p_m(K_r) = |K_r|
 \sum_{i=0}^s \frac{d_i}{|G|} \omega_i(g_u)^m
 \overline{\omega_i(g_r)}.\]
\end{lemma}

\begin{proof} If $f$ is a complex
 valued $K$ biinvariant function on $G$, define $\hat{f}(\omega_i)$
 (the Fourier transform of $f$ at the spherical function $\omega_i$)
 by \[ \hat{f}(\omega_i) = \sum_{g \in G} f(g)
 \overline{\omega_i(g)}.\] Then the Fourier inversion theorem gives
 that \[ f= \frac{1}{|G|} \sum_{i=0}^s \hat{f}(\omega_i) d_i
 \omega_i.\] It is also true that the Fourier transform of the
 convolution of two K biinvariant functions is the product of their
 Fourier transforms. Thus the m-fold convolution of $f$ is equal to \[
 \frac{1}{|G|} \sum_{i=0}^s \left(\hat{f}(\omega_i) \right)^m d_i
 \omega_i.\] The lemma now follows by taking $f=P$. \end{proof}

	To connect the study of random spherical functions with
	spectral graph theory, suppose for convenience that $(G,K)$ is
	a symmetric Gelfand pair, which means that $KgK = Kg^{-1}K$
	for all $g$. This condition holds for the examples in
	Subsections \ref{hypercube} and \ref{matching}. Then fixing a
	double coset $K_u$, one can define a graph $H_u$ whose
	vertices are the right cosets of $K$ by connecting $Kh_1$ to
	$Kh_2$ if and only if $Kh_1h_2^{-1}K = K_u$. A more general
	construction appears in \cite{Le}. The graph $H_u$ is vertex
	transitive, since $G$ acts transitively on the right cosets of
	$K$ by sending $Kh$ to $Khg^{-1}$, and $Kh_1$ is connected to
	$Kh_2$ if and only if $Kh_1g^{-1}$ is connected to
	$Kh_2g^{-1}$. Lemma \ref{specgelfand} determines the spectrum
	of random walk on $H_u$.

\begin{lemma} \label{specgelfand} Let $(G,K)$ be a symmetric Gelfand pair.
 Then for any double coset $K_u$, random walk on the graph $H_u$ has
 eigenvalues $\omega_i(g_u)$ occurring with multiplicity $d_i$ for $0
 \leq i \leq s$. \end{lemma}

\begin{proof} Let $M$ be the transition matrix for random walk on $H_u$. 
 Since $H_u$ is vertex transitive, the trace of $M^k$ is $|G|/|K|$
 multiplied by the probability that the random walk on $H_u$ started
 at the right coset $K$ is at $K$ after $k$ steps. By Theorem 7.5 of \cite{Le},
 this return probability is \[ \frac{|K|}{|G|} \sum_{i=0}^s d_i
 \omega_i(g_u)^k,\] so that the trace of $M^k$ is $\sum_{i=0}^s d_i
 \omega_i(g_u)^k$. Since this is true for all $k \geq 0$, the result
 follows. \end{proof}

	Note that since $d_0+\cdots+d_s= |G/K|$, one can define a
probability measure on $\{0,\cdots,s \}$ (or equivalently on the set
$\{\omega_0,\cdots,\omega_s\}$) by choosing $i$ with probability
$\frac{d_i |K|}{|G|}$. This is the Gelfand pair analog of the Plancherel
measure of a finite group, and in this paper it will be referred to as
Plancherel measure. Lemma \ref{specgelfand} showed that if $(G,K)$ is
a symmetric Gelfand pair, then the eigenvalues of certain graphs on
$G/K$ occur with multiplicity proportional to Plancherel measure.

\subsection{Central limit theorem for spherical functions} \label{gelfandCLT}

	The aim of this subsection is to prove a central limit theorem
	for the random variable $W$ defined by
	$W(i)=\frac{|K_u|^{1/2}}{|K|^{1/2}} \omega_i(g_u)$. Here $g_u$
	is fixed satisfying $Kg_uK= K g_u^{-1}K$ and $i$ is random
	from the Plancherel measure of the Gelfand pair $(G,K)$. It is not
	assumed that the Gelfand pair $(G,K)$ is symmetric.

	To construct an exchangeable pair to be used for Stein's
	method, it is helpful to define a Markov chain on the set of
	spherical functions of $(G,K)$. To do this fix $t$ with $1
	\leq t \leq s$ such that the spherical function $\omega_t$ is
	real-valued. Let $L_t$ be the Markov chain on
	the set $\{\omega_0,\cdots,\omega_s\}$ which transitions from
	$\omega_i$ to $\omega_j$ with probability \[
	L_t(\omega_i,\omega_j) := \frac{d_j}{|G|} \sum_{r=0}^s |K_r|
	\omega_i(g_r) \omega_t(g_r) \overline{\omega_j(g_r)}.\] As in
	the group case, one typically wants to choose $\omega_t$ as
	close as possible to the trivial spherical function.

        Lemma \ref{niceprop} verifies that $L_t$ is a Markov chain
        which is reversible with respect to Plancherel measure.

\begin{lemma} \label{niceprop} Let $t$
 with $0 \leq t \leq s$ be such that the spherical function $\omega_t$
 is real-valued. Then the transition probabilities of $L_t$ are real
 and non-negative and sum to 1. Moreover the chain $L_t$ is reversible
 with respect to the Plancherel measure of the pair
 $(G,K)$. \end{lemma}

\begin{proof}  From page 396 of \cite{Mac}, if $a_{it}^k$ are defined by \[ \omega_i
 \omega_t = \sum_{k=0}^s a_{it}^k \omega_k \] (where the notation
 $\omega_i \omega_j$ denotes the pointwise product), then $a_{it}^k$
 are real and non-negative. Hence Lemma \ref{orthog1} implies that
 $L_t(\omega_i,\omega_j)$ is real and non-negative. Since
 $\omega_j(g_0)=1$ for all $j$, Lemma \ref{orthog2} gives that \[
 \sum_{j=0}^s L_t(\omega_i,\omega_j) = \sum_{r=0}^s |K_r|
 \omega_i(g_r) \omega_t(g_r) \frac{1}{|G|} \sum_{j=0}^s d_j
 \omega_j(g_0) \overline{\omega_j(g_r)} =1.\] Reversibility of $L_t$
 with respect to Plancherel measure is equivalent to showing that \[
 \frac{d_i d_j |K|}{|G|} \sum_{r=0}^s |K_r| \omega_i(g_r) \omega_t(g_r)
 \overline{\omega_j(r)} = \frac{d_i d_j |K|}{|G|} \sum_{r=0}^s |K_r|
 \omega_j(g_r) \omega_t(g_r) \overline{\omega_i(r)}.\] Both sides are
 real by the previous paragraph, so the result follows by the
 assumption that $\omega_t$ is real valued.  \end{proof}

	An exchangeable pair $(W,W')$ to be used in a Stein's method
        approach to studying $W$ can be constructed from the chain
        $L_t$ in the usual way. First choose $i$ from Plancherel measure,
        then choose $j$ with probability $L_t(\omega_i,\omega_j)$, and
        finally let $(W,W')=(W(\omega_i),W(\omega_j))$.

\begin{lemma} \label{steinsat2} $\ee(W'|W) = \omega_t(g_u) W$.
\end{lemma}

\begin{proof} From the definitions and Lemma \ref{orthog2},
\begin{eqnarray*}
\ee(W'|\omega_i) & = & \frac{|K_u|^{1/2}}{|K|^{1/2}} \sum_{j=0}^s \frac{d_j}{|G|} \sum_{r=0}^s |K_r| \omega_i(g_r) \omega_t(g_r) \overline{\omega_j(g_r)} \omega_j(g_u)\\
& = & \frac{|K_u|^{1/2}}{|K|^{1/2}} \sum_{r=0}^s |K_r| \omega_i(g_r) \omega_t(g_r) \sum_{j=0}^s \frac{d_j}{|G|} \omega_j(g_u) \overline{\omega_j(g_r)}\\
& = &  \omega_t(g_u) W(\omega_i). \end{eqnarray*} The result follows since this depends on $\omega_i$ only through $W$. \end{proof}

	Corollary \ref{record2} is not needed in the sequel but is
	interesting.

\begin{cor}
 \label{record2} The eigenvalues of $L_t$ are $\omega_t(g_r)$ for $0
 \leq r \leq s$. The functions $\psi_r(\omega_i) =
 \frac{|K_{g_r}|^{1/2}}{|K|^{1/2}} \omega_i(g_r)$ are a basis of
 eigenvectors of $L_{\tau}$, orthonormal with respect to the inner
 product \[ \langle f_1,f_2 \rangle = \sum_{i=0}^s f_1(\omega_i)
 \overline{f_2(\omega_i)} \frac{d_i |K|}{|G|}.\]
\end{cor}

\begin{proof} By the proof of Lemma \ref{steinsat2}, $\psi_r$ is an eigenvector of $L_{t}$ with
 eigenvalue $\omega_t(g_r)$. The orthonormality assertion follows from
 Lemma \ref{orthog2}, and the basis assertion follows since the number
 of spherical functions is equal to the number of $K$ double cosets of
 $G$. \end{proof}

\begin{lemma} \label{immcor2}
$\ee(W'-W)^2 = 2(1-\omega_t(g_u))$.
\end{lemma}

\begin{proof} This is immediate from Lemmas \ref{var} and \ref{steinsat2}.
\end{proof}

\begin{lemma} \label{preGP}
$\ee((W')^2|\omega_i) = \frac{|K_u|}{|K|} \sum_{r=0}^s \omega_i(g_r) \omega_t(g_r) p_2(K_r)$.
\end{lemma}

\begin{proof} From the definitions,
\begin{eqnarray*} \ee((W')^2|\omega_i)
& = & \frac{|K_u|}{|K|} \sum_{j=0}^s \frac{d_j}{|G|} \sum_{r=0}^s |K_r| \omega_i(g_r) \omega_t(g_r) \overline{\omega_j(g_r)} \omega_j(g_u)^2\\
& = & \frac{|K_u|}{|K|} \sum_{r=0}^s |K_r| \omega_i(g_r) \omega_t(g_r) \frac{1}{|G|} \sum_{j=0}^s d_j \omega_j(g_u)^2 \overline{\omega_j(g_r)}. \end{eqnarray*} The result follows from Lemma \ref{fourier}.
\end{proof}

	Note that the sum in Lemma \ref{big1GP} begins at $r=1$ and so
	excludes the $r=0$ term which corresponds to the trivial
	spherical function.

\begin{lemma} \label{big1GP}
\[ Var(\ee[(W'-W)^2|\omega_i]) = \frac{|K_u|^2}{|K|} \sum_{r=1}^s
\frac{p_2(K_r)^2}{|K_r|} \left( \omega_t(g_r)+ 1 - 2 \omega_t(g_u) \right)^2.\]
\end{lemma}

\begin{proof} It follows from Lemmas \ref{preGP}, \ref{immcor2}, and \ref{steinsat2} that
\begin{eqnarray*}
& & Var(\ee[(W'-W)^2|\omega_i])\\
& = & \ee(\ee((W'-W)^2|\omega_i)^2) - 4(1-\omega_t(g_u))^2\\
& = & \frac{|K_u|^2}{|K|^2} \ee \left[ \sum_{r=0}^s p_2(K_r) \omega_i(g_r) \omega_t(g_r) + (1-2 \omega_t(g_u)) \omega_i(g_u)^2 \right]^2 \\
& & - 4(1-\omega_t(g_u))^2\\
& = & T_1 + T_2 + T_3 - 4(1-\omega_t(g_u))^2
\end{eqnarray*} where
\[ T_1 = \frac{|K_u|^2}{|K|^2} \ee \left( \sum_{r=0}^s p_2(K_r) \omega_i(g_r) \omega_t(g_r) \right)^2 \]
\[ T_2 = \frac{|K_u|^2}{|K|^2} 2 (1-2 \omega_t(g_u)) \sum_{r=0}^s p_2(K_r)  \omega_t(g_r) \ee \left(\omega_i(g_r)  \omega_i(g_u)^2 \right) \]
\[ T_3 = \frac{|K_u|^2}{|K|^2} (1-2 \omega_t(g_u))^2 \ee ( \omega_i(g_u)^4).\]

	For reasons similar to the finite group case, there are useful
	simplifications.  Lemma \ref{orthog2} implies that \[ T_1 =
	\frac{|K_u|^2}{|K|} \sum_{r=0}^s \frac{p_2(K_r)^2
	\omega_t(g_r)^2}{|K_r|} .\] Lemma \ref{fourier} gives that \[
	T_2 = \frac{|K_u|^2}{|K|} 2 (1-2 \omega_t(g_u)) \sum_{r=0}^s
	\frac{p_2(K_r)^2 \omega_t(g_r)}{|K_r|} .\] Also, one can write
	\[ T_3 = \frac{|K_u|^2}{|K|^2} (1-2 \omega_t(g_u))^2 p_4(K_0) =
	\frac{|K_u|^2}{|K|} (1-2 \omega_t(g_u))^2 \sum_{r=0}^s
	\frac{p_2(K_r)^2}{|K_r|}.\] The first equality used Lemma
	\ref{fourier} and the second equality used the assumption that
	$Kg_uK = Kg_u^{-1}K$. Hence \[ T_1+T_2+T_3 =
	\frac{|K_u|^2}{|K|} \sum_{r=0}^s \frac{p_2(K_r)^2}{|K_r|}
	\left( \omega_t(g_r)+1-2 \omega_t(g_u) \right)^2,\] which
	implies the result since the $r=0$ term is
	$4(1-\omega_t(g_u))^2$. \end{proof}

\begin{lemma} \label{mom1GP} Let $k$ be a positive integer.
\begin{enumerate}
\item $\ee(W'-W)^k$ is equal to \[ \left( \frac{|K_u|}{|K|} \right)^{k/2} \sum_{m=0}^k (-1)^{k-m} {k \choose m} \sum_{r=0}^s \frac{|K|}{|K_r|} \omega_t(g_r) p_m(K_r) p_{k-m}(K_r).\]
\item $\ee(W'-W)^4 = \frac{|K_u|^2}{|K|} \sum_{r=0}^s \left[ 8 \left (1-\omega_t(g_u) \right) -6 \left( 1 - \omega_t(g_r) \right) \right] \frac{p_2(K_r)^2}{|K_r|}$.
\end{enumerate}
\end{lemma}

\begin{proof} For the first assertion, note that $\ee((W'-W)^k|\omega_i)$ is equal to \begin{eqnarray*}
& & \left( \frac{|K_u|}{|K|} \right)^{k/2} \sum_{j=0}^s \frac{d_j}{|G|} \sum_{r=0}^s |K_r| \omega_i(g_r) \omega_t(g_r) \overline{\omega_j(g_r)} \left( \omega_j(g_u)-\omega_i(g_u) \right)^k\\
& = & \left( \frac{|K_u|}{|K|} \right)^{k/2} \sum_{j=0}^s \frac{d_j}{|G|} \sum_{r=0}^s |K_r| \omega_i(g_r) \omega_t(g_r) \overline{\omega_j(g_r)}\\
& & \cdot \sum_{m=0}^k (-1)^{k-m} {k \choose m} \omega_j(g_u)^m \omega_i(g_u)^{k-m}\\
& = & \left( \frac{|K_u|}{|K|} \right)^{k/2} \sum_{m=0}^k (-1)^{k-m} {k \choose m} \omega_i(g_u)^{k-m} \\
& & \cdot \sum_{r=0}^s |K_r| \omega_i(g_r) \omega_t(g_r) \sum_{j=0}^s \frac{d_j}{|G|} \omega_j(g_u)^m \overline{\omega_j(g_r)}\\
& = &  \left( \frac{|K_u|}{|K|} \right)^{k/2} \sum_{m=0}^k (-1)^{k-m} {k \choose m} \omega_i(g_u)^{k-m} \sum_{r=0}^s \omega_i(g_r) \omega_t(g_r) p_m(K_r), \end{eqnarray*} where the last equality is Lemma \ref{fourier}. Hence \begin{eqnarray*} 
\ee(W'-W)^k & = & \ee(\ee((W'-W)^k|\omega_i))\\
& = & \left( \frac{|K_u|}{|K|} \right)^{k/2} \sum_{m=0}^k (-1)^{k-m} {k \choose m}\\
& & \cdot \sum_{r=0}^s \omega_t(g_r) p_m(K_r) \sum_{i=0}^s \frac{d_i |K|}{|G|} \omega_i(g_r) \omega_i(g_u)^{k-m}. \end{eqnarray*} The first assertion follows from Lemma \ref{fourier} and the fact that $\omega_i(g_u)$ is real.

	For the second assertion, one knows from the first assertion
	that \[ \ee(W'-W)^4 = \left( \frac{|K_u|}{|K|} \right)^{2}
	\sum_{m=0}^4 (-1)^{m} {4 \choose m} \sum_{r=0}^s
	\frac{|K|}{|K_r|} \omega_t(g_r) p_m(K_r) p_{4-m}(K_r).\] The
	special case $t=0$ gives the equation \[ 0 = \left(
	\frac{|K_u|}{|K|} \right)^{2} \sum_{m=0}^4 (-1)^{m} {4 \choose
	m} \sum_{r=0}^s \frac{|K|}{|K_r|} p_m(K_r) p_{4-m}(K_r).\]
	Hence in general, $\ee(W'-W)^4$ is equal to \[ -
	\frac{|K_u|^2}{|K|^2} \sum_{m=0}^4 (-1)^{m} {4 \choose m}
	\sum_{r=0}^s \frac{|K|}{|K_r|} (1-\omega_t(g_r)) p_m(K_r)
	p_{4-m}(K_r).\] The $m=0,4$ terms both vanish since
	$p_0(K_r)=0$ if $r \neq 0$. The $m=2$ term is equal to \[ -
	\frac{6 |K_u|^2}{|K|} \sum_{r=0}^s \frac{(1-\omega_t(g_r))
	p_2(K_r)^2}{|K_r|}.\] The $m=1,3$ terms are equal and together
	contribute \begin{eqnarray*} \frac{8 |K_u|^2}{|K|}
	\sum_{r=0}^s \frac{(1-\omega_t(g_r)) p_1(K_r) p_3(K_r)}{|K_r|}
	& = & \frac{ 8 |K_u|}{|K|} (1-\omega_t(g_u)) p_3(K_u)\\ & = &
	\frac{8 |K_u|^2}{|K|^2} (1-\omega_t(g_u)) p_4(K_0)\\ & = &
	\frac{8 |K_u|^2}{|K|} (1-\omega_t(g_u)) \sum_{r=0}^s
	\frac{p_2(K_r)^2}{|K_r|}. \end{eqnarray*} This completes the
	proof of the second assertion. \end{proof}

	Arguing as in the proof of Theorem \ref{main1}, and using the
	above lemmas, one obtains the following result.

\begin{theorem}
 \label{main2} Let $(G,K)$ be a Gelfand pair, and fix a double coset
 $K_u=Kg_uK$ of $G$ satisfying $Kg_uK=Kg_u^{-1}K$. Let $\omega_t$ be a
 nontrivial spherical function which is real valued. Choosing
 $\omega_i$ from the Plancherel measure of the pair $(G,K)$, let $W=
 \frac{|K_u|^{1/2}}{|K|^{1/2}} \omega_i(g_u)$.  Then for all real
 $x_0$, \begin{eqnarray*} && \left| \pp(W \leq x_0) - \frac{1}{\sqrt{2
 \pi}} \int_{-\infty}^{x_0} e^{-\frac{x^2}{2}} dx \right|\\ & \leq &
 \frac{|K_u|}{a |K|} \sqrt{ \sum_{r=1}^s \frac{|K|p_2(K_r)^2}{|K_r|}
 \left( \omega_t(g_r)+2a-1 \right)^2 }\\ & & + \frac{1}{(\pi)^{1/4}}
 \frac{|K_u|^{1/2}}{|K|^{1/2}} \left[ \sum_{r=0}^s \left( 8 -
 \frac{6(1-\omega_t(g_r))}{a} \right) \frac{|K| p_2(K_r)^2}{|K_r|}
 \right]^{1/4}, \end{eqnarray*} where
 $a=1-\omega_t(g_u)$. \end{theorem}

\subsection{Example: The Hypercube} \label{hypercube}

	The $n$ dimensional hypercube $\zz_2^n$ consists of $n$-tuples
        of 0's and 1's. Random walk on it proceeds by picking a random
        coordinate and changing it. The spectrum of this random walk
        is well known; for each $0 \leq i \leq n$ there is an
        eigenvalue $1-\frac{2i}{n}$ occurring with multiplicity ${n
        \choose i}$ (see for instance page 28 of \cite{Di}). Thus by
        the usual central limit theorem for the binomial distribution, the
        spectrum of the hypercube is asymptotically normal with an
        error term $O(n^{-1/2})$. The purpose of this subsection is to
        revisit this classical result from the viewpoint of Gelfand
        pairs, illustrating the construction of Subsection
        \ref{gelfandCLT}.

	To begin, note from the remarks on page 58 of \cite{Di}, that
        the hypercube can be viewed as $G/K$ for a certain Gelfand
        pair $(G,K)$. Namely $G$ is the semidirect product of $\zz_2^n$
        with $S_n$, where the group multiplication is
        $(x,\pi)(y,\tau)=(x+\pi(y),\pi \tau)$, where $\tau(x)$
        permutes the coordinates of $x$. $K$ is the subgroup
        $\{(0,\pi):\pi \in S_n \}$. The induced module $1_K^G$
        decomposes as $\oplus_{r=0}^n V_r$ where $d_r={n \choose
        r}$. Thus the Plancherel measure chooses $i \in
        \{0,\cdots,n\}$ with probability $\frac{{n \choose
        i}}{2^n}$. For $0 \leq r \leq n$, the double coset $K_r$ in
        $G$ consists of elements $(x,\pi)$ where $x$ has r coordinates
        equal to 1. Thus $|K_r| = {n \choose r} n!$, and as usual let
        $g_r$ denote some element of $K_r$. The spherical
        function $\omega_i$ ($0 \leq i \leq n$) is given by \[
        \omega_i(g_r) = \frac{1}{{n \choose i}} \sum_{m=0}^i (-1)^m {r
        \choose m} {n-r \choose i-m},\] which is a Krawtchouk
        polynomial if one overlooks the ${n \choose i}$ in the
        denominator.

	The following result emerges from Theorem \ref{main2}. Since
	$\omega_i(g_1)=1-\frac{2i}{n}$ (or by Lemma
	\ref{specgelfand}), it can be seen as a central limit theorem
	for the spectrum of random walk on the hypercube.

\begin{theorem} \label{hypbound1} Let $W=\sqrt{n}
 \omega_i(g_1)$ where $i$ is chosen from Plancherel measure. Then for
 all real $x_0$, \[ \left| \pp(W \leq x_0) - \frac{1}{\sqrt{2 \pi}}
 \int_{-\infty}^{x_0} e^{-\frac{x^2}{2}} dx \right| \leq \left(
 \frac{8}{\pi n} \right)^{1/4} .\] \end{theorem}

\begin{proof} Apply
 Theorem \ref{main2} with $u=1$ and $t=1$. Then $a=\frac{2}{n}$. One
 computes that $p_2(K_0)=\frac{1}{n}$, $p_2(K_2)=1-\frac{1}{n}$, and
 that $p_2(K_r)=0$ for all $r \neq 0,2$. Since
 $\omega_1(g_2)=1-\frac{4}{n}$, one has that
 $\omega_1(g_2)+2a-1=0$. Thus the first error term in Theorem
 \ref{main2} is 0. The second error term in Theorem \ref{main2} is
 computed to be $\left( \frac{8}{\pi n} \right)^{1/4}$, implying the
 result. \end{proof}

	To get $O(n^{-1/2})$ bounds by Stein's method, it will be
	shown that $|W'-W|$ is bounded, so that one can use the
	version of Stein's method in Theorem \ref{rinrot} instead of
	that in Theorem \ref{steinbound}. The boundedness of $|W'-W|$
	will follow from Lemma \ref{birth}, which proves that $L_1$ is
	in fact a birth-death chain.

\begin{lemma} \label{birth} The chain $L_1$ on the set $\{0,\cdots,n\}$ is a birth-death chain with transition probabilities \[ L_1(i,j) = \left\{ \begin{array}{ll} \frac{i}{n} & \mbox{if $j=i-1$}\\
1-\frac{i}{n} & \mbox{if $j=i+1$} \end{array} \right.\] \end{lemma}

\begin{proof} From the three term recurrence for Krawtchouk polynomials on page 152 of \cite{MS}, it follows that
\[ (i+1){n \choose i+1} \omega_{i+1}(g_r) = (n-2r) {n \choose i}  \omega_i(g_r)- (n-i+1) {n \choose i-1}  \omega_{i-1}(g_r).\] Since $n-2r = n \omega_1(g_r)$, one obtains that \[ n {n \choose i} \omega_1(g_r) \omega_i(g_r) = (i+1) {n \choose
 i+1} \omega_{i+1}(g_r) + (n-i+1) {n \choose i-1} \omega_{i-1}(g_r).\]
 Simplifying one obtains the relation
\[ \omega_1(g_r) \omega_i(g_r) = \left( 1- \frac{i}{n} \right) \omega_{i+1}(g_r) + \left( \frac{i}{n} \right) \omega_{i-1}(g_r).\] The result now follows from Lemma \ref{orthog1} and the definition of $L_1$. \end{proof}

	Now an $O(n^{-1/2})$ error term is established.

\begin{theorem}
 \label{hypbound2} Let $W=\sqrt{n} \omega_i(g_1)$ where $i$ is chosen
 from Plancherel measure. Then for
 all real $x_0$, \[ \left| \pp(W \leq x_0) - \frac{1}{\sqrt{2 \pi}}
 \int_{-\infty}^{x_0} e^{-\frac{x^2}{2}} dx \right| \leq 5 n^{-1/2}
 .\] \end{theorem}

\begin{proof} Apply the variation of Theorem \ref{main2} which would arise from using Theorem \ref{rinrot}
 instead of Theorem \ref{steinbound}. Recall that
 $a=\frac{2}{n}$. Also $|W'-W| \leq A$ with $A=\frac{2}{\sqrt{n}}$
 since $L_1$ is a birth-death chain. As explained in the proof of
 Theorem \ref{hypbound1}, $Var(\ee[(W'-W)^2|W])=0$. The result
 follows. \end{proof}

	A a final remark, Lemma \ref{birth} shows that $L_1$ is
closely related to the birth-death chain used by Stein \cite{Stn} in
proving a central limit theorem for $X_1+...+X_n$, where the $X's$ are
independent and each equal to 0 or 1 with probability 1/2. To form an
exchangeable pair Stein chose a random $l \in \{1,\cdots,n\}$ and then
replaced $X_l$ by a new random variable which is also 0 or 1 with
probability 1/2 and independent of all of the $X's$. The chain $L_1$
would arise by picking a random index $l$ and switching the value of
$X_l$.

\subsection{Example: Graphs on Perfect Matchings} \label{matching}

	In this example $G=S_{2n}$ and $K$ is the hyperoctahedral
	group of signed permutations on $n$ symbols, of size $2^n
	n!$. This Gelfand pair is discussed at length in Chapter 7 of
	\cite{Mac} and in Section 3 of \cite{HSS}, to which we refer
	the reader for proofs of facts in the next two
	paragraphs. Another useful reference is \cite{DHol}.

	The induced representation $1_K^G$ decomposes as
	$\oplus_{|\lambda|=n} V_{\lambda}$, where $\lambda$ ranges
	over all partitions of size $n$. An explicit formula for the
	numbers $d_{\lambda}$ appears in \cite{Mac} and is not needed
	in what follows. However it is worth remarking that the Plancherel measure of
	$(G,K)$, which chooses $\lambda$ with probability $\frac{2^n
	n!  d_{\lambda}}{(2n)!}$, is the $\alpha=2$ case of the so
	called Jack$_{\alpha}$ measure on partitions, which chooses
	$\lambda$ with probability \[ \frac{\alpha^n n!}{\prod_{s \in
	\lambda} (\alpha a(s) + l(s)+1) (\alpha a(s) + l(s) +
	\alpha)} \] where the product is over all boxes of the
	partition. Here $a(s)$ is the number of boxes in the same row
	of $s$ and to the right of $s$ (the ``arm'' of $s$), and
	$l(s)$ is the number of boxes in the same column of $s$ and
	below $s$ (the ``leg'' of $s$). For example the partition of 5
	below 
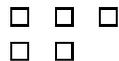
\[
\begin{array}{c c c} \framebox{}& \framebox{}& \framebox{} \\
\framebox{}& \framebox{}& \end{array} \] would have Jack$_{\alpha}$
measure \[ \frac{60 \alpha^2}{(2 \alpha+2)(3 \alpha+1) (\alpha+2)(2
\alpha+1)(\alpha+1)}.\] The Jack measure on partitions is of
interest to researchers in random matrix theory \cite{BO1},\cite{O2},
\cite{K2}.

	The double cosets $K_{\mu}$ of $K$ in $G$ are parameterized by
	partitions $\mu$ of size $n$ and have the following concrete
	description. A perfect matching of $\{1,\cdots,2n\}$ can be
	regarded as a 1-regular graph with vertex set
	$\{1,\cdots,2n\}$. Let $\epsilon$ be the ``identity matching''
	in which $i$ is adjacent to $n+i$ for $1 \leq i \leq n$. Given
	a permutation $w$ in $S_{2n}$, let $\delta(w)$ be the perfect
	matching of $\{1,\cdots,2n\}$ in which $i$ is adjacent to $j$
	if and only $|w(i)-w(j)|=n$. For example,
	$\delta(id)=\epsilon$. Note that the union $\delta_1 \cup
	\delta_2$ of two 1-regular graphs is a 2-regular graph, and
	thus a disjoint union of even length cycles. Let
	$\Lambda(\delta_1,\delta_2)$ be the partition whose parts are
	half the cycle lengths of $\delta_1 \cup \delta_2$. Then
	$Kw_1K = Kw_2 K$ if and only if $\Lambda(\epsilon,\delta(w_1))
	= \Lambda(\epsilon,\delta(w_2))$. Thus $Kw^{-1}K = KwK$ for
	all $w$, so that the machinery of Subsection \ref{gelfandCLT}
	is applicable. One also has that $Kw_1=Kw_2$ if and only if
	$\delta(w_1)=\delta(w_2)$. It follows that $G/K$ is in
	bijection with the perfect matchings of $\{1,\cdots,2n\}$ and
	that $\frac{|K_{\mu}|}{|K|}$ is equal to the number of perfect
	matchings $\delta$ with $\Lambda(\epsilon,\delta)=\mu$, which
	by elementary counting is $\frac{2^n n!}{2^{l(\mu)} \prod_j
	m_j(\mu)! j^{m_j(\mu)}}$ where $l(\mu)$ is the number of parts
	of $\mu$ and $m_j(\mu)$ is the number of parts of $\mu$ of size $j$.

	The main purpose of this subsection is to prove a central
	limit theorem for the random variable $W = \sqrt{2^{i-1} {n
	\choose i} (i-1)!}  \omega_{\lambda}(g_{(i,1^{n-i})})$, where
	$i$ is fixed and $\lambda$ is chosen from the Plancherel
	measure of $(G,K)$. Then Lemma \ref{specgelfand} gives a
	central limit theorem for the spectrum of random walk on the
	graph $H_{(i,1^{n-i})}$, whose vertices are the perfect
	matchings of $\{1,\cdots,2n\}$, with an edge between matchings
	$\delta_1$ and $\delta_2$ if and only if
	$\Lambda(\delta_1,\delta_2) = (i,1^{n-i})$. In the case $i=2$,
	the spectrum of this graph was determined in \cite{DHol}, and
	was shown to be asymptotically normal in \cite{F2} (with error
	term $O(n^{-1/4})$) and then in \cite{CF}, \cite{F3} (with
	error term $O(n^{-1/2})$). But the arguments in those papers
	used different information, and it was not clear that
	they could be pushed through to larger $i$. Theorem \ref{main2}
	will be used to deduce a central limit theorem for any fixed
	$i$ with error term $O(n^{-1/4})$.

	To apply Theorem \ref{main2}, $\omega_t$ will be taken to be
	the spherical function $\omega_{(n-1,1)}$. An explicit formula
	for $\omega_{(n-1,1)}$ is available.

	\begin{lemma} \label{valspher} (\cite{Mac}, p. 411) Let
	$m_1(\mu)$ denote the number of parts of size 1 of $\mu$. Then
	\[ \omega_{(n-1,1)}(g_{\mu}) =
	\frac{(2n-1)m_1(\mu)-n}{2n(n-1)} \] for all $\mu$ of size
	$n$. \end{lemma}

	Lemma \ref{con} is helpful.

	\begin{lemma} \label{con} (\cite{HSS}, Lemma 3.2) The
	coefficient of $K_{\mu}$ in $K_{\tau} K_{(i,1^{n-i})}$ is
	equal to $\frac{|K|^2}{|K_{\mu}|}$ multiplied by the number of
	pairs of perfect matchings $(\delta,\gamma)$ such that
	$\Lambda(\epsilon,\delta) = (i,1^{n-i})$,
	$\Lambda(\delta,\gamma)=\tau$, and
	$\Lambda(\epsilon,\gamma) = \mu$.  \end{lemma}

	The final combinatorial ingredient is an analog of Lemma \ref{term1}.

	\begin{lemma} \label{term1GP} Consider the random walk on $G$
	generated by $K_{(i,1^{n-i})}$. Then for $i$ fixed and $n \geq
	2i$, $\frac{p_2(K_{\mu})^2 |K|}{|K_{\mu}|}$ is equal to
	\begin{enumerate}
\item $\frac{i^2}{4^{i-1}} n^{-2i} + O(n^{-2i-1})$ if $K_{\mu}=K_{(1^n)}=K$.
\item $\frac{2i^2}{4^{i-1}} n^{-2i} + O(n^{-2i-1})$ if $K_{\mu}=K_{(i,i,1^{n-2i})}$.
\item $O(n^{-2i-1})$ otherwise.
\end{enumerate} 
\end{lemma}

\begin{proof} Lemma \ref{con} implies that $\frac{p_2(K_{\mu}) |K_{(i,1^{n-i})}|^2}{|K|^2}$
is equal to the number of pairs of matchings $(\delta,\gamma)$ such
that $\Lambda(\epsilon,\delta) = (i,1^{n-i})$,
$\Lambda(\delta,\gamma)=(i,1^{n-i})$, and $\Lambda(\epsilon,\gamma) = \mu$.

For the first assertion, it follows either from the previous paragraph
or from Lemmas \ref{fourier} and \ref{orthog2} that $p_2(K) =
\frac{1}{2^{i-1} {n \choose i} (i-1)!}$. For the second assertion, it
is straightforward from the previous paragraph that
$p_2(K_{(i,i,1^{n-2i})}) = 1 + O(n^{-1})$. Then use the formula for
$\frac{|K_{\mu}|}{|K|}$ given earlier in this subsection.

For the third assertion, there are two cases. The first case is that
$\mu$ has $n-2i$ parts of size 1, but is not equal to
$(i,i,1^{n-2i})$. Then $p_2(K_{\mu})=O(n^{-1})$ and the formula for
$\frac{|K|}{|K_{\mu}|}$ shows it to be at most $c_i n^{-2i}$ where
$c_i$ is a constant depending on $i$. So in the first case, the result
is proved. The second case is that $\mu$ has $n-2i+r$ parts of size 1,
where $1 \leq r < 2i$. From the first paragraph of the proof, it is
straightforward to see that $p_2(K_{\mu})=O(n^{-r})$. Also the formula
for $\frac{|K|}{|K_{\mu}|}$ shows it to be $O(n^{r-2i})$, proving the
result. \end{proof}

	Combining the ingredients, one deduces the following result.

\begin{theorem} \label{CLTgel} Choose $\lambda$ from the Plancherel
 measure of the Gelfand pair $(G,K)$, and define a random variable $W=
 \sqrt{2^{i-1} {n \choose i} (i-1)!}
 \omega_{\lambda}(g_{(i,1^{n-i})})$. Then there is a constant $A_i$
 such that for all real $x_0$, \[ \left| \pp(W \leq x_0) -
 \frac{1}{\sqrt{2 \pi}} \int_{-\infty}^{x_0} e^{-\frac{x^2}{2}} dx
 \right| \leq A_i n^{-1/4} .\] \end{theorem}

\begin{proof} One applies Theorem \ref{main2}, choosing $u=(i,1^{n-i})$ and $\omega_t$ to be the spherical
 function $\omega_{(n-1,1)}$. By Lemma \ref{valspher}, $\omega_t$ is
 real valued and $a=\frac{i(2n-1)}{2n(n-1)}$. Observe that
 $\omega_t(g_{\mu})+2a-1$ vanishes if $\mu=(i,i,1^{n-2i})$ and by
 Lemma \ref{valspher} is $O(n^{-1})$ for any other $\mu$ such that
 $p_2(K_{\mu}) \neq 0$, since any such $\mu$ has at least $n-2i$ parts
 of size 1. Thus by Lemma \ref{term1GP}, \[ \sum_{\mu \neq (1^n)} \frac{|K|
 p_2(K_{\mu})^2}{|K_{\mu}|} (\omega_t(g_{\mu}) + 2a -1)^2 =
 O(n^{-2i-3}).\] Since $\frac{|K_u|}{|K|} = 2^{i-1} {n \choose i}
 (i-1)!$, the first error term in Theorem \ref{main2} is at most $A_i'
 n^{-1/2}$, where $A_i'$ is a constant depending on $i$.

	 To bound the second error term in Theorem \ref{main2}, one
 applies Lemma \ref{term1GP} to see that \[ \sum_{\mu}
 \left(8 - \frac{6(1-\omega_t(g_{\mu}))}{a} \right) \frac{|K|
 p_2(K_{\mu})^2}{|K_{\mu}|} = O(n^{-2i-1}).\] Note that cancellation
 occurred for the coefficient of $n^{-2i}$ coming from the terms
 $\mu=(1^n)$ and $\mu=(i,i,1^{n-2i})$. It follows that the second
 error term is at most $A_i'' n^{-1/4}$, where $A_i''$ is another
 constant depending only on $i$. \end{proof}

\section{Twisted Gelfand Pairs} \label{twisted}

	Let $G$ be a finite group and $H$ a subgroup of $G$. If
	$\phi$ is a linear character of $H$ such that the induced
	representation $Ind_H^G(\phi)$ is multiplicity free, then
	$(G,H,\phi)$ is referred to as a twisted Gelfand pair (this
	terminology was introduced in \cite{Stem}).

	The approach taken in this section is a bit different than
	that in the sections on finite groups and Gelfand pairs. Due
	both to a lack of interesting examples which are not Gelfand
	pairs and to technical complications which do not arise for
	finite groups or Gelfand pairs, a completely general theory is
	not developed. Instead, we focus on one very interesting
	example: character values of random projective representations
	of the symmetric group. It should however be noted that many
	of the lemmas resemble those of earlier sections, and the
	calculations are organized in a way which should generalize to
	other examples.

	It is known (see \cite{HH} for a friendly exposition) that the
	character values in question are expressed in terms of the
	coefficients $X^{\lambda}_{\rho}$, where $\lambda$ is a strict
	(i.e. all parts distinct) partition of $n$ that parameterizes
	an irreducible character, and $\mu$ is an odd (i.e. all parts
	odd) partition of $n$ that parameterizes a conjugacy
	class. Throughout this section $DP(n)$ will denote the set of
	strict partitions of size $n$, and $OP(n)$ will denote the set
	of odd partitions of size $n$.

	Central limit theorems for random projective representations
	of the symmetric group have been studied by Ivanov
	\cite{I}. The underlying probability measure on strict
	partitions chooses $\lambda$ with probability \[
	\frac{2^{n-l(\lambda)} g_{\lambda}^2}{n!}, \] where
	$l(\lambda)$ is the number of parts of a partition $\lambda$
	and $g_{\lambda} = X^{\lambda}_{(1^n)}$ is the number of
	standard shifted tableaux of shape $\lambda$
	(\cite{HH},\cite{Mac}). This measure on strict partitions is
	known as shifted Plancherel measure, and Ivanov proves the
	following result.

\begin{theorem}
 \label{projivan} (\cite{I}) Fix $i \geq 1$. Let $\lambda$ be chosen
 from the shifted Plancherel measure on strict partitions of size $n$. Then
 as $n \rightarrow \infty$, the random variable
 $\frac{n^{\frac{2i+1}{2}} X^{\lambda}_{(2i+1,1^{n-2i-1})}}{2^i \sqrt{2i+1}
 g_{\lambda}}$ converges in distribution to a normal random variable
 with mean 0 and variance 1.
\end{theorem}

	This section refines the result of Ivanov, which was proved by
	the method of moments, so as to obtain an error
	term. In the statement of the result, recall that $z_{\mu} =
	\prod_{j \geq 1} j^{m_j(\mu)} m_j(\mu)!$, where $m_j(\mu)$ is
	the number of parts of $\mu$ of size $j$.

\begin{theorem} \label{projerror} Fix $i \geq 1$ and
 let $\mu=(2i+1,1^{n-2i-1})$. Choosing $\lambda$ from the shifted
 Plancherel measure on partitions of size $n$, define a random
 variable \[ W = \sqrt{\frac{n!}{z_{\mu} 2^{n-l(\mu)}}}
 \frac{X^{\lambda}_{\mu}}{g_{\lambda}}.\] Then there is a constant
 $A_i$ such that for all real $x_0$, \[ \left| \pp(W \leq x_0) -
 \frac{1}{\sqrt{2 \pi}} \int_{-\infty}^{x_0} e^{-\frac{x^2}{2}} dx
 \right| \leq A_i n^{-1/4}.\]
\end{theorem}

	The organization of this section is as follows. Subsection
	\ref{trtgelfand} collects and develops background from
	representation theory. Subsection \ref{exchange} defines and
	studies a Markov chain to be used in the construction of an
	exchangeable pair for a Stein's method proof of Theorem
	\ref{projerror}. This is more subtle than the corresponding
	treatment in earlier sections and involves interesting
	combinatorics, since the ``obvious'' adaptation of the
	construction for Gelfand pairs does not work. Subsection
	\ref{twCLT} studies the exchangeable pair arising from the
	Markov chain in Section \ref{exchange}, and uses it to prove
	Theorem \ref{projerror}.

\subsection{Background from representation theory} \label{trtgelfand}

	If $G$ is a finite group, $K$ a subgroup of $G$, and $\phi$ a
	linear character of $K$ such that $Ind_K^G(\phi)$ is
	multiplicity free, the triple $(G,K,\phi)$ is called a twisted
	Gelfand pair. The Hecke algebra of the triple $(G,K,\phi)$ is
	the $\cc G$ subalgebra $e \cc G e$, where $e$ is the primitive
	idempotent of $\cc K$ defined by \[ e =\frac{1}{|K|} \sum_{k
	\in K} \phi(k^{-1}) k.\] One reason that twisted Gelfand pairs
	are interesting is that the Hecke algebra of the triple
	$(G,K,\phi)$ is commutative. The paper \cite{Stem} is a good
	reference for the theory of twisted Gelfand pairs.

	The twisted Gelfand pair of interest to us is the one studied
	in \cite{Stem}. Thus $G=S_{2n}$ and $K$ is the hyperoctahedral
	group $B_n$ of signed permutations, imbedded in $G$ as the
	centralizer of the involution $(1,2)(3,4) \cdots
	(2n-1,2n)$. To define $\phi$, note that $B_n$ is the
	semidirect product of the groups $T_n$ and $\Sigma_n$, where
	$T_n$ is the subgroup (isomorphic to $\zz_2^n$) generated by
	$(1,2),\cdots,(2n-1,2n)$ and $\Sigma_n$ is the subgroup
	(isomorphic to $S_n$) generated by the ``double
	transpositions'' $(2i-1,2j-1)(2i,2j)$ for $1 \leq i<j \leq
	n$. Then $\phi$ is the linear character of $B_n$ whose
	restriction to $\Sigma_n$ is the sign character, and whose
	restriction to $T_n$ is trivial.

	Stembridge \cite{Stem} defines twisted spherical functions for
	twisted Gelfand pairs (the analog of spherical functions for
	Gelfand pairs). One of his main results is that for
	$(S_{2n},B_n,\phi)$, their values are (aside from scalar
	multiples) equal to the $X^{\lambda}_{\mu}$, where $\lambda$ is
	a strict partition and $\mu$ is an odd partition.

	Next we record orthogonality relations for the quantities
	$X^{\lambda}_{\mu}$. These orthogonality relations are
	a special case of more general orthogonality relations for
	coefficients of power sum symmetric functions in
	Hall-Littlewood polynomials, where the parameter in the Hall
	Littlewood polynomial is $-1$. What this paper calls
	$X^{\lambda}_{\mu}$ is written as $X^{\lambda}_{\mu}(-1)$ in
	Chapter 3 of \cite{Mac}.

\begin{lemma} \label{torthog1} (\cite{Mac}, page 247) For $\lambda,\rho \in DP(n)$, \[ \sum_{\mu \in OP(n)} \frac{2^{l(\mu)}}{z_{\mu}} X^{\lambda}_{\mu} X^{\rho}_{\mu} = \delta_{\lambda,\rho} 2^{l(\lambda)}.\]
\end{lemma}

\begin{lemma} \label{torthog2} (\cite{Mac}, page 247) For $\mu,\sigma \in OP(n)$, \[ \sum_{\lambda \in DP(n)} \frac{1}{2^{l(\lambda)}} X^{\lambda}_{\mu} X^{\lambda}_{\sigma} = \delta_{\mu,\sigma} \frac{z_{\mu}}{2^{l(\mu)}}.\]
\end{lemma}

	As in Subsection \ref{matching}, the double cosets of $B_n$ in
	$S_{2n}$ are indexed by partitions $\nu$ of $n$. It is
	useful to specify representatives $w_{\nu}$. For the case
	$\nu=(n)$, one defines \[ w_{(n)} = (1,2,\cdots,2n),\] and for
	the general case $\nu=(\nu_1,\cdots,\nu_l)$, one defines \[
	\omega_{\nu} = w_{(\nu_1)} \circ \cdots \circ w_{(\nu_l)} \]
	where the operation $x \circ y$ (for $x \in S_{2i}, y \in
	S_{2j}$) denotes the embedding of $S_{2i} \times S_{2j}$ in
	$S_{2i+2j}$ with $S_{2i}$ acting on $\{1,\cdots,2i\}$ and
	$S_{2j}$ acting on $\{2i+1,\cdots,2i+2j\}$.

	For $\nu \in OP(n)$, define \[ \tilde{K}_{\nu} =
	\frac{1}{|B_n|^2} \sum_{x_1,x_2 \in B_n} \phi(x_1x_2) x_1
	w_{\nu} x_2.\] Corollary 3.2 of \cite{Stem} shows that if $\nu
	\not \in OP(n)$, then $\tilde{K}_{\nu}=0$, and that
	$\{\tilde{K}_{\nu}: \nu \in OP(n)\}$ is a basis for the Hecke
	algebra of $(S_{2n},B_{n},\phi)$. Thus for $\mu, \nu \in OP(n)$, it is
	natural to study the coefficient of $\tilde{K}_{\nu}$ in
	$(\tilde{K}_{\mu})^m$. Lemma \ref{tcountsol} gives
	a character theoretic expression for this coefficient and is
	analogous to Lemmas \ref{countsol} and \ref{fourier}.

\begin{lemma}
 \label{tcountsol} Suppose that $\mu_1,\cdots,\mu_m, \nu \in
OP(n)$. Then the coefficient of $\tilde{K}_{\nu}$ in $\tilde{K}_{\mu_1} \cdots
\tilde{K}_{\mu_m}$ is equal to \[ \frac{2^{l(\mu_1)-n} \cdots
2^{l(\mu_m)-n}}{z_{\nu}} \sum_{\rho \in DP(n)} 2^{n-l(\rho)}
g_{\rho}^2 \frac{X^{\rho}_{\nu}}{g_{\rho}}
\frac{X^{\rho}_{\mu_1}}{g_{\rho}} \cdots
\frac{X^{\rho}_{\mu_m}}{g_{\rho}} .\] In particular, the coefficient
of $\tilde{K_{\nu}}$ in $(\tilde{K}_{\mu})^m$ is \[ \frac{n!
(2^{l(\mu)-n})^m}{z_{\nu}} \ee\left[ \frac{X^{\rho}_{\nu}}{g_{\rho}}
\left( \frac{X^{\rho}_{\mu}}{g_{\rho}} \right)^m \right] , \] where
$\rho$ is random from shifted Plancherel measure on partitions of size
$n$. \end{lemma}

\begin{proof}
 The equality is immediate from the definition of shifted
 Plancherel measure, so it is enough to establish the first
 expression. Stembridge \cite{Stem} provides a basis of orthogonal
 idempotents $\{ E_{\rho}: \rho\in DP(n) \}$ of the Hecke algebra of
 $(G,H,\phi)$ in terms of the projective characters of the
 symmetric group. More precisely, it follows from Proposition 4.1 and
 Corollary 6.2 of \cite{Stem} that one obtains such a basis of
 orthogonal idempotents by defining
\[ E_{\rho} = 2^{n-l(\rho)} g_{\rho} \sum_{\mu \in OP(n)} \frac{X^{\rho}_{\mu}}{z_{\mu}} \tilde{K}_{\mu}. \]

	Multiplying both sides of this equation by $2^{l(\mu)-n}
\frac{X^{\rho}_{\mu}}{g_{\rho}}$ and summing over all $\rho \in
DP(n)$, it follows from Lemma \ref{torthog2} that \[ \tilde{K}_{\mu} =
2^{l(\mu)-n} \sum_{\rho \in DP(n)} \frac{X^{\rho}_{\mu}}{g_{\rho}}
E_{\rho}.\] Thus \[ \tilde{K}_{\mu_1} \cdots \tilde{K}_{\mu_m} = 2^{l(\mu_1)-n} \cdots
2^{l(\mu_m)-n} \sum_{\rho \in DP(n)} \frac{X^{\rho}_{\mu_1}}{g_{\rho}}
\cdots \frac{X^{\rho}_{\mu_m}}{g_{\rho}} E_{\rho} .\] The result now
follows by the formula for $E_{\rho}$ in the previous
paragraph. \end{proof}

	Lemma \ref{valtwist} derives an explicit formula for
	$X^{(n-1,1)}_{\mu}$. This is crucial to our approach.

\begin{lemma} \label{valtwist} Let $m_1(\mu)$ denote the
 number of parts of size 1 of $\mu$. Suppose that $n \geq 3$, so that
 $(n-1,1)$ is a strict partition. Then \[ X^{(n-1,1)}_{\mu} = m_1(\mu)
 - 2 \] for all odd partitions $\mu$ of size $n$.
\end{lemma}

\begin{proof}  Recall the definition of power sum symmetric functions:
 for $i \geq 1$, one sets $p_i = \sum_j x_j^i$, and for $\mu$ a
 partition of $n$, one sets $p_{\mu} = \prod_i p_i^{m_i(\mu)}$. The
 argument also uses Schur's Q-functions (\cite{Mac}, Sec. 3.8). More
 precisely, equation 7.5 on page 247 of \cite{Mac} shows that for $\mu
 \in OP(n)$ and $\lambda \in DP(n)$, the coefficient of the power sum
 symmetric function $p_{\mu}$ in $Q_{\lambda}$ is equal to
 $\frac{2^{l(\mu)} X^{\lambda}_{\mu}}{z_{\mu}}$.

	 To proceed one needs an expression for $Q_{(n-1,1)}$. Using
 the notation that $[x^n] f(x)$ is the coefficient of $x^n$ in $f(x)$,
 one has by page 253 of \cite{Mac} that \begin{eqnarray*} Q_{(n-1,1)}
 & = & [t_1^{n-1} t_2] \left(\frac{1-t_2/t_1}{1+t_2/t_1} \right)
 \prod_{i=1}^2 \prod_{j=1}^{\infty} \frac{1+t_ix_j}{1-t_ix_j}\\ & = & [t_1^{n-1}]
 \prod_{j=1}^{\infty} \frac{1+t_1x_j}{1-t_1x_j} [t_2] \prod_{j=1}^{\infty}
 \frac{1+t_2x_j}{1-t_2x_j} -2 [t_1^n] \prod_{j=1}^{\infty}
 \frac{1+t_1x_j}{1-t_1x_j}\\ & = & 2p_1 Q_{n-1} - 2
 Q_{n}. \end{eqnarray*}

	By page 248 of \cite{Mac}, $X^{n}_{\mu}=1$ for all $\mu \in
	OP(n)$. It follows from the first paragraph that for $\mu \in
	OP(n)$, the coefficient of $p_{\mu}$ in $Q_n$ is
	$\frac{2^{l(\mu)}}{z_{\mu}}$ and (by considering separately
	the cases that $m_1(\mu)=0$ and $m_1(\mu) \neq 0$) that the
	coefficient of $p_{\mu}$ in $p_1 Q_{n-1}$ is $\frac{m_1(\mu)
	2^{l(\mu)-1}}{z_{\mu}}$. Thus the coefficient of $p_{\mu}$ in
	$Q_{(n-1,1)}$ is $\frac{(m_1(\mu)-2) 2^{l(\mu)}}{z_{\mu}}$,
	which by the first paragraph proves the result. \end{proof}

\subsection{Markov chains on strict partitions}
 \label{exchange}

	This subsection discusses a Markov chain on $DP(n)$ to be used
	in defining an exchangeable pair for the proof of Theorem
	\ref{projerror}. 

	Motivated by constructions in the group and Gelfand pair
	cases, it would be natural, for $\tau$ a fixed element of
	$DP(n)$, to define a ``Markov chain'' $J_{\tau}$ on $DP(n)$
	with transition ``probabilities'' \[ J_{\tau}(\lambda,\rho) :=
	\frac{g_{\rho}}{2^{l(\rho)} g_{\lambda} g_{\tau}} \sum_{\nu
	\in OP(n)} \frac{2^{l(\nu)} X^{\lambda}_{\nu} X^{\rho}_{\nu}
	X^{\tau}_{\nu}}{z_{\nu}} .\] Using Lemma \ref{torthog2}, one
	can see that $\sum_{\rho} J_{\tau}(\lambda,\rho) = 1$ for all
	$\lambda$, and that $J_{\tau}$ satisfies the reversibility
	condition \[ \frac{2^{n-l(\lambda)} g_{\lambda}^2}{n!} J_{\tau}(\lambda,\rho) = \frac{2^{n-l(\rho)} g_{\rho}^2}{n!}
	J_{\tau}(\rho,\lambda) \] for all $\lambda,\rho \in
	DP(n)$. Unfortunately, the quantity $J_{\tau}(\lambda,\rho)$
	can be negative. For instance one can check from Lemma
	\ref{valtwist} that $J_{(2,1)}((2,1),(2,1))=-1$.

	To deal with the complication raised in the previous
	paragraph, it is helpful to introduce a genuine Markov chain
	$L$ on $DP(n)$, with transition probabilities \[ L(\lambda,
	\rho) = \frac{2g_{\rho}}{n g_{\lambda}} \sum_{\eta \in DP(n-1)
	\atop \eta \nearrow \lambda, \eta \nearrow \rho} 2^{l(\eta)-l(\rho)}.\] Here
	$\eta \nearrow \lambda$ means that $\eta$ is obtained from
	$\lambda$ by decreasing the size of some part by exactly one.

	Lemma \ref{ischain} proves that $L$ is a Markov chain which is
	reversible with respect to shifted Plancherel measure. As is
	mentioned in the proof, the definition of $L$ was motivated by
	the theory of harmonic functions on Bratteli diagrams.

\begin{lemma} \label{ischain} The transition
 probabilities of $L$ are real and non-negative and sum to 1. Moreover
 the chain $L$ is reversible with respect to shifted Plancherel
 measure.
\end{lemma}

\begin{proof}
 This is a special case of a construction in Section 2 of
 \cite{F1}. To see this, one takes the underlying Bratteli diagram to
 be the Schur graph, the properties of which are discussed in Section 5 of
 \cite{BO2}. The vertices of the Schur graph are all partitions of all
 non-negative integers with distinct parts, and the edge multiplicity
 between $\eta \in DP(n-1)$ and $\lambda \in DP(n)$ is 1 if $\eta
 \nearrow \lambda$, and is 0 otherwise. The combinatorial dimension of
 a shape $\lambda$ is $g_{\lambda}$, and if $\pi$ denotes shifted
 Plancherel measure, the function $\frac{\pi(\lambda)}{g_{\lambda}}$
 is harmonic on the Schur graph. \end{proof}

	Proposition \ref{2chains} will establish a fundamental
	relation between the Markov chain $L$ and the ``Markov chain''
	$J_{(n-1,1)}$. The argument involves symmetric function
	theory, and first a lemma is needed about properties of a
	certain subring $\Gamma$ of the ring of symmetric functions,
	defined on page 252 of \cite{Mac}.

	Recall from page 255 of \cite{Mac} that there is an inner
	product on $\Gamma$ which satisfies the properties \begin{enumerate} \item $\langle p_{\lambda},p_{\mu} \rangle = 2^{-l(\lambda)} z_{\lambda}
	\delta_{\lambda,\mu}$ if $\lambda,\mu \in OP(n)$.
\item $\langle P_{\lambda},P_{\mu} \rangle= 2^{-l(\lambda)}
\delta_{\lambda,\mu}$ if $\lambda,\mu \in DP(n)$. \end{enumerate} Here
$p_{\lambda}$ is a power sum symmetric function and $P_{\lambda}$ is a
Hall-Littlewood polynomial with the parameter $t=-1$. It is also
helpful to recall, from page 247 of \cite{Mac}, that for $\lambda
\in DP(n)$, \[ P_{\lambda} = 2^{-l(\lambda)} \sum_{\nu \in OP(n)}
\frac{2^{l(\nu)}}{z_{\nu}} X^{\lambda}_{\nu} p_{\nu}.\]

\begin{lemma} \label{symlem} Let $p_1^{\perp}$ be the adjoint in $\Gamma$ of multiplication by $p_1$.
\begin{enumerate}
\item (\cite{Mac}, page 265) $p_1^{\perp} = \frac{1}{2}
\frac{\partial}{\partial p_1}$.
\item $p_1^{\perp} P_{\lambda} = \sum_{\eta \in DP(n-1) \atop \eta
\nearrow \lambda} 2^{l(\eta)-l(\lambda)} P_{\eta}$.
\end{enumerate}
\end{lemma}

\begin{proof} Only
 the second assertion needs to be proved. Consider the coefficient of
 $P_{\eta}$ in $p_1^{\perp} P_{\lambda}$.  It is \[ \frac{ \langle p_1^{\perp}
 P_{\lambda}, P_{\eta} \rangle}{\langle P_{\eta},P_{\eta} \rangle} = \frac{\langle  P_{\lambda},
 p_1 P_{\eta} \rangle}{\langle P_{\eta},P_{\eta} \rangle} = 2^{l(\eta)} \langle P_{\lambda}, p_1
 P_{\eta} \rangle.\] From Section 3.8 of \cite{Mac}, \[ p_1 P_{\eta} =
 \sum_{\lambda \in DP(n) \atop \eta \nearrow \lambda} P_{\lambda},\]
 which implies the result. \end{proof}

	Proposition \ref{2chains} can now be proved.

\begin{prop} \label{2chains} Suppose that $n \geq 3$, so that $(n-1,1)$ is a strict partition. Then \[ L(\lambda,\rho) = \frac{n-2}{n} J_{(n-1,1)}(\lambda,\rho)\] for $\lambda,\rho \in DP(n)$ such that $\lambda \neq \rho$. 
\end{prop}

\begin{proof} By Lemma \ref{valtwist}, \[ J_{(n-1,1)}(\lambda,\rho) = 	\frac{g_{\rho}}{2^{l(\rho)} g_{\lambda} (n-2)} \sum_{\nu
	\in OP(n)} \frac{2^{l(\nu)} X^{\lambda}_{\nu} X^{\rho}_{\nu}
	(m_1(\nu)-2)}{z_{\nu}}.\] Since $\lambda \neq \rho$, Lemma
	\ref{torthog1} and part 1 of Lemma \ref{symlem} imply that
	this is equal to \begin{eqnarray*} & &
	\frac{g_{\rho}}{2^{l(\rho)} g_{\lambda} (n-2)} \sum_{\nu \in
	OP(n)} \frac{2^{l(\nu)} X^{\lambda}_{\nu} X^{\rho}_{\nu}
	m_1(\nu)}{z_{\nu}}\\ & = & \frac{g_{\rho}}{2^{l(\rho)}
	g_{\lambda} (n-2)} \left \langle \sum_{\nu \in OP(n)} m_1(\nu)
	\frac{2^{l(\nu)}}{z_{\nu}} X^{\lambda}_{\nu} p_{\nu},
	\sum_{\nu \in OP(n)} \frac{2^{l(\nu)}}{z_{\nu}} X^{\rho}_{\nu}
	p_{\nu} \right \rangle\\ & = & \frac{2^{l(\lambda)}
	g_{\rho}}{g_{\lambda} (n-2)} \langle 2 p_1 p_1^{\perp}
	P_{\lambda},P_{\rho} \rangle\\ & = & \frac{2^{l(\lambda)+1}
	g_{\rho}}{g_{\lambda} (n-2)} \langle p_1^{\perp} P_{\lambda},
	p_1^{\perp} P_{\rho} \rangle. \end{eqnarray*}

	By part 2 of Lemma \ref{symlem}, this is \begin{eqnarray*} & &
	\frac{2^{l(\lambda)+1} g_{\rho}}{g_{\lambda} (n-2)} \left
	\langle \sum_{\eta \in DP(n-1) \atop \eta \nearrow \lambda}
	2^{l(\eta)-l(\lambda)} P_{\eta} , \sum_{\eta \in DP(n-1) \atop
	\eta \nearrow \rho} 2^{l(\eta)-l(\rho)} P_{\eta} \right
	\rangle \\ & = & \frac{2g_{\rho}}{g_{\lambda} (n-2)}
	\sum_{\eta \in DP(n-1) \atop \eta \nearrow \lambda, \eta
	\nearrow \rho} 2^{l(\eta)-l(\rho)} \\ & = & \frac{n}{n-2}
	L(\lambda,\rho). \end{eqnarray*}
\end{proof}

\subsection{Central limit theorem for shifted Plancherel measure} \label{twCLT}

	This subsection studies the statistic $W=\sqrt{
	\frac{n!}{z_{\mu} 2^{n-l(\mu)}}}
	\frac{X^{\lambda}_{\mu}}{g_{\lambda}}$, where $\mu \in OP(n)$
	is fixed and $\lambda \in DP(n)$ is chosen from shifted
	Plancherel measure. The main goal is a proof of Theorem
	\ref{projerror}.

	Using $L$, one constructs an exchangeable pair $(W,W')$ as
	follows. Choose $\lambda$ from the shifted Plancherel measure on
	partitions of size $n$. Then choose $\rho$ with probability
	$L(\lambda,\rho)$ and let $(W,W')=(W(\lambda),W(\rho))$. Using
	$J_{\tau}$, one constructs an ``exchangeable pair'' $(W,W^*)$
	as follows. Choose $\lambda$ from the shifted Plancherel measure
	on partitions of size $n$. Then choose $\rho$ with ``probability''
	$J_{\tau}(\lambda,\rho)$ and let
	$(W,W^*)=(W(\lambda),W(\rho))$. The pair $(W,W')$ is a valid
	candidate for Stein's method, but the pair $(W,W^*)$ is much
	easier to work with, and by Proposition \ref{2chains}, when
	$\tau=(n-1,1)$, this gives insight into the genuine
	exchangeable pair $(W,W')$. Even though the transition
	probabilities of $J_{\tau}$ can be negative, for convenience
	the usual language of probability theory (expected value,
	variance, etc.) will be used when working with them.

\begin{lemma} \label{tsteinsat2}
\begin{enumerate}
\item $\ee(W^*|W) = \left( \frac{X^{\tau}_{\mu}}{g_{\tau}} \right)W$.
\item $\ee(W'|W) = \frac{m_1(\mu)}{n} W$.
\end{enumerate}
\end{lemma}

\begin{proof} For the first assertion, the definition of $J_{(n-1,1)}$ and Lemma \ref{torthog2} imply
that \begin{eqnarray*} \ee(W^*|\lambda) & = & \sqrt{\frac{n!}{z_{\mu}
2^{n-l(\mu)}}} \sum_{\rho \in DP(n)} \frac{g_{\rho}}{2^{l(\rho)}
g_{\lambda} g_{\tau}} \sum_{\nu \in OP(n)} \frac{2^{l(\nu)}
X^{\lambda}_{\nu} X^{\rho}_{\nu} X^{\tau}_{\nu}}{z_{\nu}}
\frac{X^{\rho}_{\mu}}{g_{\rho}}\\ & = & \sqrt{\frac{n!}{z_{\mu}
2^{n-l(\mu)}}} \frac{1}{g_{\lambda} g_{\tau}} \sum_{\nu \in OP(n)}
\frac{2^{l(\nu)} X^{\lambda}_{\nu} X^{\tau}_{\nu}}{z_{\nu}} \sum_{\rho
\in DP(n)} \frac{ X^{\rho}_{\nu} X^{\rho}_{\mu} }{2^{l(\rho)}} \\ & =
& \sqrt{\frac{n!}{z_{\mu} 2^{n-l(\mu)}}} \frac{X^{\lambda}_{\mu}}{g_{\lambda}}
\frac{X^{\tau}_{\mu}}{g_{\tau}}\\ & = & \left( \frac{X^{\tau}_{\mu}}{g_{\tau}}
\right) W(\lambda). \end{eqnarray*} The first assertion follows since this depends on $\lambda$ only through $W$.

	For the second assertion, observe that by Proposition
\ref{2chains}, \begin{eqnarray*} \ee(W'-W|\lambda) & = & \sum_{\rho} L(\lambda,\rho) \left(W(\rho)-W(\lambda) \right)\\ & = &
\frac{n-2}{n} \sum_{\rho} J_{(n-1,1)}(\lambda,\rho)
\left( W(\rho) - W(\lambda) \right)\\ & = & - \left( \frac{n-2}{n}
\right) W(\lambda) + \frac{n-2}{n} \sum_{\rho} J_{(n-1,1)}(\lambda,\rho)
W(\rho). \end{eqnarray*} By the first assertion and Lemma
\ref{valtwist}, this is \[ -\left(\frac{n-2}{n} \right) W(\lambda) +
\frac{n-2}{n} \frac{m_1(\mu)-2}{n-2} W(\lambda) = \left( \frac{m_1(\mu)}{n} -1
\right) W(\lambda), \] which implies the result.  \end{proof}

	Corollary \ref{diago} will not be needed but is worth recording.

\begin{cor} \label{diago} The eigenvalues of $J_{\tau}$ are $\frac{X^{\tau}_{\mu}}{g_{\tau}}$ as $\mu$ ranges over $OP(n)$. The functions $\psi_{\mu}(\lambda) = \sqrt{\frac{n!}{z_{\mu} 2^{n-l(\mu)}}} \frac{X^{\lambda}_{\mu}}{g_{\lambda}}$ are a basis of eigenvectors of $J_{\tau}$, orthonormal with respect to the inner product \[ \langle f_1,f_2 \rangle = \sum_{\lambda \in DP(n)} f_1(\lambda) \overline{f_2(\lambda)} \frac{2^{n-l(\lambda)} g_{\lambda}^2}{n!}.\]
\end{cor}

\begin{proof} The proof of part 1 of Lemma \ref{tsteinsat2} shows that $\psi_{\mu}$ is an eigenvector of $J_{\tau}$ with eigenvalue $\frac{X^{\tau}_{\mu}}{g_{\tau}}$. The orthonormality assertion follows from Lemma \ref{torthog2} and the fact from \cite{Mac} that all $X^{\lambda}_{\mu}$ are real valued. The basis assertion follows since $|DP(n)|=|OP(n)|$. \end{proof}

\begin{lemma} \label{timmcor2} $\ee(W^*-W)^2 = 2 \left(1- \frac{X^{\tau}_{\mu}}{g_{\tau}} \right)$. \end{lemma}

\begin{proof} This follows from the proof of Lemma \ref{var} (which does not require non-negative transition probabilities) and part 1 of Lemma \ref{tsteinsat2}. 
\end{proof}

	For the remainder of this subsection, $p_m(\tilde{K}_{\nu})$
	will denote the coefficient of $\tilde{K}_{\nu}$ in
	$(\tilde{K}_{\mu})^m$. When $m=0$ this is to be interpreted
	through Lemma \ref{tcountsol}, so that
	$p_0(\tilde{K}_{(1^n)})=1$ and $p_0(\tilde{K}_{\nu}) = 0$ for
	$\nu \neq (1^n)$. Due to the signs in the definition of
	$\tilde{K}_{\nu}$, these numbers are not probabilities so care
	must be taken in working with them. For instance it is not
	true that $\sum_{\nu \in OP(n)} p_m(\tilde{K}_{\nu})
	=1$. However the following three relations will be useful.

\begin{lemma} \label{sumneq1} $\sum_{\nu \in OP(n)}  p_2(\tilde{K}_{\nu}) 2^{l(\nu)-n} = 2^{2[l(\mu)-n]}$. \end{lemma}

\begin{proof} Using Lemma \ref{tcountsol},
 one has that, \[ \sum_{\nu \in OP(n)} p_2(\tilde{K}_{\nu})
 2^{l(\nu)} = 2^{2l(\mu)-n} \sum_{\rho \in DP(n)} 2^{-l(\rho)}
 \frac{(X^{\rho}_{\mu})^2}{g_{\rho}} \sum_{\nu \in OP(n)}
 \frac{2^{l(\nu)}}{z_{\nu}} X^{\rho}_{\nu} .\] By page 248 of
 \cite{Mac}, $X^{n}_{\nu}=1$, so that Lemma \ref{torthog1} implies
 that \[ \sum_{\nu \in OP(n)} \frac{2^{l(\nu)}}{z_{\nu}} X^{\rho}_{\nu}
 X^{n}_{\nu} = 2 \delta_{\rho,n}.\] The result
 follows. \end{proof}

\begin{lemma} \label{4step} $\ee \left( \frac{X^{\lambda}_{\mu}}{g_{\lambda}} \right)^4 = (2^{n-l(\mu)})^4 \sum_{\nu \in OP(n)} p_2(\tilde{K}_{\nu})^2 \frac{2^{l(\nu)} z_{\nu}}{2^n n!}$.
 \end{lemma}

\begin{proof} Consider
 the coefficient of $\tilde{K}_{(1^n)}$ in $(\tilde{K}_{\mu})^4$. On
 one hand, by Lemma \ref{tcountsol} it is equal to $(2^{l(\mu)-n})^4
 \ee \left(\frac{X^{\lambda}_{\mu}}{g_{\lambda}} \right)^4$. On the
 other hand, \[ (\tilde{K}_{\mu})^4 = (\tilde{K}_{\mu})^2
 (\tilde{K}_{\mu})^2 = \left[ \sum_{\nu \in OP(n)} p_2(\tilde{K}_{\nu})
 \tilde{K}_{\nu} \right]^2.\] From Lemmas \ref{tcountsol} and Lemma
 \ref{torthog2}, it follows that the coefficient of
 $\tilde{K}_{(1^n)}$ in $\tilde{K}_{\nu_1} \tilde{K}_{\nu_2}$ is 0 if
 $\nu_1 \neq \nu_2$, and that the coefficient of $\tilde{K}_{(1^n)}$ in
 $(\tilde{K}_{\nu})^2$ is $\frac{2^{l(\nu)-n} z_{\nu}}{n!}$. Thus the
 coefficient of $\tilde{K}_{(1^n)}$ in $(\tilde{K}_{\mu})^4$ is equal
 to \[ \sum_{\nu \in OP(n)} p_2(\tilde{K}_{\nu})^2 \frac{2^{l(\nu)-n}
 z_{\nu}}{n!}.\] Comparing the two expressions for the coefficient of
 $\tilde{K}_{(1^n)}$ in $(\tilde{K}_{\mu})^4$ proves the
 result. \end{proof}

\begin{lemma} \label{3step} $p_3(\tilde{K}_{\mu}) = \frac{2^{n-l(\mu)}} {z_{\mu}} \sum_{\nu \in OP(n)} p_2(\tilde{K}_{\nu})^2 2^{l(\nu)-n} z_{\nu}$.
\end{lemma}

\begin{proof} By Lemma \ref{tcountsol}, \[ p_3(\tilde{K}_{\mu})=\frac{n! (2^{l(\mu)-n})^3}{z_{\mu}} \ee \left( \frac{X^{\lambda}_{\mu}}{g_{\lambda}}
 \right)^4.\] Now use Lemma \ref{4step}. \end{proof}

	The next lemmas are crucial.

\begin{lemma} \label{tpreGP} $\ee((W^*)^2|\lambda) = \frac{n!}{z_{\mu} 2^{l(\mu)} } \sum_{\nu \in OP(n)} \frac{X^{\lambda}_{\nu}}{g_{\lambda}} \frac{X^{\tau}_{\nu}}{g_{\tau}} 2^{l(\nu)} p_2(\tilde{K}_{\nu})$.
\end{lemma}

\begin{proof} It follows from the definitions that $
\ee((W^*)^2|\lambda)$ is equal to
\begin{eqnarray*} &  & \frac{n!}{z_{\mu} 2^{n-l(\mu)}} \sum_{\rho
\in DP(n)} \frac{g_{\rho}}{2^{l(\rho)} g_{\lambda} g_{\tau}} \sum_{\nu
\in OP(n)} \frac{2^{l(\nu)} X^{\lambda}_{\nu} X^{\rho}_{\nu}
X^{\tau}_{\nu}}{z_{\nu}} \left( \frac{X^{\rho}_{\mu}}{g_{\rho}} \right)^2\\ & = &
\frac{n!}{z_{\mu} 2^{n-l(\mu)}} \sum_{\nu \in OP(n)}
\frac{X^{\lambda}_{\nu}}{g_{\lambda}}
\frac{X^{\lambda}_{\tau}}{g_{\tau}} 2^{l(\nu)} \left(
\frac{1}{z_{\nu}} \sum_{\rho \in DP(n)} \frac{g_{\rho}^2}{2^{l(\rho)}}
\frac{X^{\rho}_{\nu}}{g_{\rho}} \left( \frac{X^{\rho}_{\mu}}{g_{\rho}}
\right)^2 \right). \end{eqnarray*} The result now follows from Lemma
\ref{tcountsol}.
\end{proof}

\begin{lemma}
 \label{tbig1GP} \[ Var(\ee[(W^*-W)^2|\lambda]) = \frac{n!
 2^n}{z_{\mu}^2 2^{2 l(\mu)}} \sum_{\nu \in OP(n) \atop \nu \neq
 (1^n)} 2^{l(\nu)} p_2(\tilde{K}_{\nu})^2 z_{\nu} \left(
 \frac{X^{\tau}_{\nu}}{g_{\tau}}+1-\frac{2 X^{\tau}_{\mu}}{g_{\tau}}
 \right)^2\]
\end{lemma}

\begin{proof} By Lemmas \ref{tpreGP}, \ref{timmcor2}, and \ref{tsteinsat2}, $Var(\ee[(W^*-W)^2|\lambda])$ is equal to
\begin{eqnarray*}
 &  & \ee(\ee((W^*-W)^2|\lambda)^2) - 4
\left(1-\frac{X^{\tau}_{\mu}}{g_{\tau}} \right)^2\\ & = & \ee \left[
\frac{n!}{z_{\mu} 2^{l(\mu)}} \sum_{\nu \in OP(n)}
\frac{X^{\lambda}_{\nu}}{g_{\lambda}} \frac{X^{\tau}_{\nu}}{g_{\tau}}
2^{l(\nu)} p_2(\tilde{K}_{\nu}) + \left(1 - \frac{2
X^{\tau}_{\mu}}{g_{\tau}} \right) W^2 \right]^2\\ & & - 4
\left(1-\frac{X^{\tau}_{\mu}}{g_{\tau}} \right)^2\\ & = & T_1 + T_2 +
T_3 - 4 \left(1-\frac{X^{\tau}_{\mu}}{g_{\tau}} \right)^2\\
\end{eqnarray*} where

\[ T_1 = \left( \frac{n!}{z_{\mu} 2^{l(\mu)}} \right)^2 \ee \left( \sum_{\nu \in OP(n)} p_2(\tilde{K}_{\nu}) \frac{X^{\lambda}_{\nu}}{g_{\lambda}} \frac{X^{\tau}_{\nu}}{g_{\tau}} 2^{l(\nu)} \right)^2 \]

\[ T_2 = 2 \left( 1-\frac{2 X^{\tau}_{\mu}}{g_{\tau}} \right) \frac{(n!)^2}{z_{\mu}^2} \sum_{\nu \in OP(n)} \frac{X^{\tau}_{\nu} p_2(\tilde{K}_{\nu})}{g_{\tau} 2^{n-l(\nu)}} \ee \left( \frac{X^{\lambda}_{\nu} (X^{\lambda}_{\mu})^2}{(g_{\lambda})^3} \right) \]

\[ T_3 = \left( 1-\frac{2 X^{\tau}_{\mu}}{g_{\tau}} \right)^2 \left( \frac{n!}{2^{n-l(\mu)} z_{\mu}} \right)^2 \ee \left( \frac{X^{\lambda}_{\mu}}{g_{\lambda}} \right)^4. \]

	These expressions for $T_1,T_2,T_3$ can be simplified. Lemma
	\ref{torthog2} implies that \[ T_1 = \frac{n!
	2^n}{(z_{\mu})^2 2^{2 l(\mu)}} \sum_{\nu \in OP(n)} p_2(\tilde{K}_{\nu})^2
	2^{l(\nu)} z_{\nu} \left( \frac{X^{\tau}_{\nu}}{g_{\tau}}
	\right)^2.\] From Lemma \ref{tcountsol} it follows that \[ T_2
	= 2 \left( 1-\frac{2 X^{\tau}_{\mu}}{g_{\tau}} \right) \frac{n!
	2^n}{(z_{\mu})^2 2^{2 l(\mu)}} \sum_{\nu \in OP(n)} p_2(\tilde{K}_{\nu})^2
	2^{l(\nu)} z_{\nu} \frac{X^{\tau}_{\nu}}{g_{\tau}}.\] By
	Corollary \ref{4step}, \[ T_3 = \frac{n! 2^n}{2^{2 l(\mu)}
	(z_{\mu})^2} \left( 1-\frac{2 X^{\tau}_{\mu}}{g_{\tau}} \right)^2
	\sum_{\nu \in OP(n)} p_2(\tilde{K}_{\nu})^2 2^{l(\nu)} z_{\nu}.\] Hence, one can
	write \[ T_1+T_2+T_3 = \frac{n! 2^n}{(z_{\mu})^2 2^{2 l(\mu)}}
	\sum_{\nu \in OP(n)} 2^{l(\nu)} z_{\nu} p_2(\tilde{K}_{\nu})^2 \left(
	\frac{X^{\tau}_{\nu}}{g_{\tau}} + 1 - 
	\frac{2 X^{\tau}_{\mu}}{g_{\tau}} \right)^2 \] and the result
	follows since by Lemmas \ref{tcountsol} and \ref{torthog2},
	the $\nu=(1^n)$ term in the above sum is $4 \left(1 -
	\frac{X^{\tau}_{\mu}}{g_{\tau}} \right)^2$. \end{proof}

	Lemma \ref{tmom1GP} will also be helpful.

\begin{lemma} \label{tmom1GP} Let $k$ be a positive integer.
\begin{enumerate}
\item $\ee(W^*-W)^k$ is equal to \[ \left( \frac{n!
2^{n-l(\mu)}}{z_{\mu}} \right)^{k/2}  \sum_{m=0}^k
(-1)^{k-m} {k \choose m} \sum_{\nu \in OP(n)} \frac{X^{\tau}_{\nu}
z_{\nu}}{g_{\tau} n! 2^{n-l(\nu)}} p_m(\tilde{K}_{\nu}) p_{k-m}(\tilde{K}_{\nu}).\]
\item $ \ee(W^*-W)^4$ is equal to \[ \left( \frac{n!
2^{n-l(\mu)}}{z_{\mu}} \right)^2 \left[ \sum_{\nu \in OP(n)} \left( 8
\left( 1-\frac{X^{\tau}_{\mu}}{g_{\tau}} \right)-6 \left(
1-\frac{X^{\tau}_{\nu}}{g_{\tau}} \right) \right) \frac{z_{\nu}
p_2(\tilde{K}_{\nu})^2}{n! 2^{n-l(\nu)}} \right].\]
\end{enumerate}
\end{lemma}

\begin{proof} To prove the first assertion, observe that
\begin{eqnarray*}
& & \ee((W^*-W)^k|\lambda)\\ & = & \left(
\frac{n!}{z_{\mu} 2^{n-l(\mu)}} \right)^{k/2} \sum_{\rho \in DP(n)}
\frac{g_{\rho}}{2^{l(\rho)} g_{\lambda} g_{\tau}} \sum_{\nu \in OP(n)}
\frac{2^{l(\nu)} X^{\lambda}_{\nu} X^{\rho}_{\nu}
X^{\tau}_{\nu}}{z_{\nu}}\\ & & \cdot \sum_{m=0}^k (-1)^{k-m} {k
\choose m} \left(\frac{X^{\rho}_{\mu}}{g_{\rho}} \right)^m
\left(\frac{X^{\lambda}_{\mu}}{g_{\lambda}} \right)^{k-m}\\ & = &
\left( \frac{n!}{z_{\mu} 2^{n-l(\mu)}} \right)^{k/2}
\frac{1}{g_{\lambda} g_{\tau}} \sum_{m=0}^k (-1)^{k-m} {k \choose m}
\left(\frac{X^{\lambda}_{\mu}}{g_{\lambda}} \right)^{k-m}\\ & & \cdot
\sum_{\nu \in OP(n)} \frac{2^{l(\nu)} X^{\lambda}_{\nu}
X^{\tau}_{\nu}}{z_{\nu}} \sum_{\rho \in DP(n)}
\frac{g_{\rho}^2}{2^{l(\rho)}} \frac{X^{\rho}_{\nu}}{g_{\rho}} \left(
\frac{X^{\rho}_{\mu}}{g_{\rho}} \right)^m \\ & = & \left(
\frac{n!}{z_{\mu} 2^{n-l(\mu)}} \right)^{k/2} \frac{1}{2^n g_{\lambda}
g_{\tau}} \sum_{m=0}^k (-1)^{k-m} {k
\choose m}  (2^{n-l(\mu)})^m \left(\frac{X^{\lambda}_{\mu}}{g_{\lambda}} \right)^{k-m}\\
& & \cdot \sum_{\nu \in OP(n)} 2^{l(\nu)} X^{\lambda}_{\nu}
X^{\tau}_{\nu} p_m(\tilde{K}_{\nu}),
\end{eqnarray*} where the last equation is Lemma \ref{tcountsol}. Hence $\ee(W^*-W)^k$ is equal to 
\begin{eqnarray*}
& & \ee(\ee((W^*-W)^k|\lambda))\\
 & = & \left( \frac{n!}{z_{\mu}
2^{n-l(\mu)}} \right)^{k/2} \frac{(2^{n-l(\mu)})^k}{2^n} \sum_{m=0}^k
(-1)^{k-m} {k \choose m}\\ & & \cdot \sum_{\nu \in OP(n)} 2^{l(\nu)} \frac{
X^{\tau}_{\nu}}{g_{\tau}} \frac{ p_m(\tilde{K}_{\nu}) }{(2^{n-l(\mu)})^{k-m}} \ee \left[
\frac{X^{\lambda}_{\nu}}{g_{\lambda}} \left(
\frac{X^{\lambda}_{\mu}}{g_{\lambda}} \right)^{k-m}
\right]. \end{eqnarray*} The first assertion now follows from Lemma
\ref{tcountsol}.

	For the second assertion, the first assertion gives
 that $\ee(W^*-W)^4$ is equal to \[ \left(\frac{n!
 2^{n-l(\mu)}}{z_{\mu}} \right)^2 \sum_{m=0}^4 (-1)^{m} {4 \choose m}
 \sum_{\nu \in OP(n)} \frac{z_{\nu}}{2^{n-l(\nu)} n!}
 \frac{X^{\tau}_{\nu}}{g_{\tau}} p_m(\tilde{K}_{\nu})
 p_{4-m}(\tilde{K}_{\nu}) .\] The special case $\tau=(n)$ gives that
 \[ 0 = \left(\frac{n!  2^{n-l(\mu)}}{z_{\mu}} \right)^2 \sum_{m=0}^4
 (-1)^{m} {4 \choose m} \sum_{\nu \in OP(n)}
 \frac{z_{\nu}}{2^{n-l(\nu)} n!}  p_m(\tilde{K}_{\nu})
 p_{4-m}(\tilde{K}_{\nu}).\] Thus in general, $\ee(W^*-W)^4$ is equal
 to \[ - \left(\frac{n!  2^{n-l(\mu)}}{z_{\mu}} \right)^2 \sum_{m=0}^4
 (-1)^{m} {4 \choose m} \sum_{\nu \in OP(n)}
 \left(1-\frac{X^{\tau}_{\nu}}{g_{\tau}} \right) \frac{z_{\nu}
 p_m(\tilde{K}_{\nu}) p_{4-m}(\tilde{K}_{\nu})}{2^{n-l(\nu)} n!}  .\]
 The $m=0,4$ terms vanish since only $\nu=(1^n)$ could contribute, but
 it contributes 0. The $m=2$ term contributes \[ -6 \left(\frac{n!
 2^{n-l(\mu)} }{z_{\mu}} \right)^2 \sum_{\nu \in OP(n)} \left(
 1-\frac{X^{\tau}_{\nu}}{g_{\tau}} \right) \frac{z_{\nu}
 p_2(\tilde{K}_{\nu})^2}{n!  2^{n-l(\nu)}}.\] The $m=1,3$ terms are
 equal and together contribute \[ 8 \left(\frac{n!
 2^{n-l(\mu)}}{z_{\mu}} \right) \left(
 1-\frac{X^{\tau}_{\mu}}{g_{\tau}} \right) p_3(\tilde{K}_{\mu}).\] By
 Lemma \ref{3step}, this is \[ 8 \left( \frac{n!
 2^{n-l(\mu)}}{z_{\mu}} \right)^2 \left( 1 -
 \frac{X^{\tau}_{\mu}}{g_{\tau}} \right) \sum_{\nu \in OP(n)}
 \frac{z_{\nu} p_2(\tilde{K}_{\nu})^2}{n! 2^{n-l(\nu)}} .\] Adding the
 terms together proves the second assertion.  \end{proof}

	The final combinatorial ingredient for the proof of Theorem
	\ref{projerror} is an upper bound for
	$p_2(\tilde{K}_{\nu})$. Lemma \ref{reduce} reduces this to a
	result obtained earlier in the section on Gelfand pairs.

\begin{lemma}
 \label{reduce} Let $p_2(\tilde{K}_{\nu})$ be as in this subsection and
 let $p_2(K_{\nu})$ be as in Subsection \ref{matching}. Then
 $p_2(\tilde{K}_{\nu}) \leq p_2(K_{\nu})$ for all $\nu \in OP(n)$.
\end{lemma}

\begin{proof}
 Clearly $p_2(K_{\nu})$ is equal to $|B_n \omega_{\nu} B_n|$
 multiplied by the coefficient of $\omega_{\nu}$ in \[ \left(
 \frac{1}{|B_n|^2} \sum_{x_1,x_2 \in B_n} x_1 \omega_{\mu} x_2
 \right)^2.\] Next consider $p_2(\tilde{K}_{\nu})$. By definition it
 is the coefficient of $\tilde{K}_{\nu}$ in
 $(\tilde{K}_{\mu})^2$. Corollary 3.2 of \cite{Stem} gives that the
 coefficient of $\omega_{\nu}$ in $\tilde{K}_{\nu}$ is $\frac{1}{|B_n
 \omega_{\nu} B_n|}$. Hence $p_2(\tilde{K}_{\nu})$ is $|B_n
 \omega_{\nu} B_n|$ multiplied by the coefficient of $\omega_{\nu}$ in
 \[ \left( \frac{1}{|B_n|^2} \sum_{x_1,x_2 \in B_n} \phi(x_1 x_2) x_1
 \omega_{\mu} x_2 \right)^2,\] where $\phi$ is the linear character of
 $B_n$ described in Subsection \ref{rtgelfand}. Since $\phi(x_1
 x_2)$ is always $\pm 1$, the result follows. \end{proof}

	Finally, we prove Theorem \ref{projerror}.

\begin{proof}
 (Of Theorem \ref{projerror}) One can assume that $n \geq 2(2i+1)$. In
 the $O$ notation throughout the proof, $i$ is fixed and $n$ is
 growing.

	 One applies Theorem \ref{steinbound} to the pair $(W,W')$,
 where $\mu=(2i+1,1^{n-2i-1})$. Then $a=1-\frac{m_1(\mu)}{n}$ by part
 2 of Lemma \ref{tsteinsat2}. Proposition \ref{2chains} implies that
 \[ \ee((W'-W)^2|\lambda) = \frac{n-2}{n} \ee((W^*-W)^2|\lambda)\] for
 all $\lambda$. Thus \[ Var(\ee[(W'-W)^2|\lambda]) = \left(
 \frac{n-2}{n} \right)^2 Var(\ee[(W^*-W)^2|\lambda]).\] Lemmas
 \ref{tbig1GP} and \ref{valtwist} give that the first error term in
 applying Theorem \ref{steinbound} to the pair $(W,W')$ is at most \[
 \frac{n-2}{n} \frac{n! 2^{n-l(\mu)}}{\left(
 1-\frac{m_1(\mu)}{n} \right) z_{\mu}} \sqrt{ \sum_{\nu \in OP(n) \atop \nu \neq (1^n)}
 \frac{z_{\nu} p_2(\tilde{K}_{\nu})^2}{n!  2^{n-l(\nu)}} \left(
 \frac{n+m_1(\nu)-2m_1(\mu)}{n-2} \right)^2}.\] Observe that
 $\frac{n-2}{n} \frac{n! 2^{n-l(\mu)}}{\left(
 1-\frac{m_1(\mu)}{n} \right) z_{\mu}}$ is $O(n^{2i+2})$. In the sum
 over $\nu \neq (1^n)$, the term coming from
 $\nu=(2i+1,2i+1,1^{n-4i-2})$ contributes 0. Since $n \geq 2(2i+1)$,
 Lemmas \ref{reduce} and \ref{term1GP} imply that $\frac{z_{\nu}
 p_2(\tilde{K}_{\nu})^2}{n! 2^{n-l(\nu)}} = O(n^{-4i-3})$ for all
 other $\nu$ in the sum, and the number of such $\nu$ is bounded by a
 constant depending on $i$ but not $n$. Moreover if
 $p_2(\widetilde{K}_{\nu}) \neq 0$, then $p_2(K_{\nu}) \neq 0$,
 implying by Lemma \ref{con} that $m_1(\nu) \geq n-4i-2$ and thus that
 $ \frac{n+m_1(\nu)-2m_1(\mu)}{n-2}=O(n^{-1})$. Combining these
 observations gives that the first term in the upper bound of Theorem
 \ref{steinbound} is $O(n^{-1/2})$.

	By the Cauchy-Schwarz inequality, \[ \ee|W'-W|^3 \leq \sqrt{
	\ee(W'-W)^2 \ee(W'-W)^4} .\] Proposition \ref{2chains} implies
	that $\ee(W'-W)^k = \frac{n-2}{n} \ee(W^*-W)^k$ for all
	$k$. Hence Lemmas \ref{timmcor2} and \ref{valtwist} imply that
	the second error term in Theorem \ref{steinbound} is at most
	\[ \left[ \frac{(n-2)}{\pi n}
	\frac{\ee(W^*-W)^4}{(1-\frac{m_1(\mu)}{n})} \right]^{1/4} .\]
	By part 2 of Lemma \ref{tmom1GP} and Lemma \ref{valtwist},
	this can be written as \[ \left[ \frac{1}{\pi} \frac{(n!
	2^{n-l(\mu)})^2}{(1-\frac{m_1(\mu)}{n}) (z_{\mu})^2} \sum_{\nu
	\in OP(n)} \left( 2 - \frac{8m_1(\mu)}{n} + \frac{6
	m_1(\nu)}{n} \right) \frac{z_{\nu} p_2(\tilde{K}_{\nu})^2}{n!
	2^{n-l(\nu)}} \right]^{1/4}. \] The quantity $\frac{(n!
	2^{n-l(\mu)})^2}{(1-\frac{m_1(\mu)}{n}) (z_{\mu})^2}$ is
	$O(n^{4i+3})$. The contribution from $\nu=(1^n)$ to the sum is
	$\frac{8(2i+1)}{n} p_2(\tilde{K}_{(1^n)})^2$, which by Lemmas
	\ref{tcountsol} and \ref{torthog2} is $\frac{8
	(2i+1)^3}{2^{4i}} n^{-4i-3} + O(n^{-4i-4})$. By Lemmas
	\ref{reduce} and \ref{term1GP}, the contribution to the sum
	from $\nu$ not equal to either $(1^n)$ or
	$(2i+1,2i+1,1^{n-4i-2})$ is $O(n^{-4i-4})$. The contribution
	from $\nu=(2i+1,2i+1,1^{n-4i-2})$ is $-\frac{4(2i+1)}{n}
	\frac{z_{\nu} p_2(\tilde{K}_{\nu})^2}{n! 2^{n-l(\nu)}}$. Lemma
	\ref{sumneq1} implies that
	$p_2(\tilde{K}_{(2i+1,2i+1,1^{n-4i-2})})= 1+O(n^{-1})$, since
	$p_2(\tilde{K}_{\nu}) = O(n^{-1})$ for $\nu \neq
	(2i+1,2i+1,1^{n-4i-2})$ by Lemma \ref{reduce}. Thus the
	contribution from $\nu=(2i+1,2i+1,1^{n-4i-2})$ is
	$-\frac{8(2i+1)^3}{2^{4i}} n^{-4i-3} + O(n^{-4i-4})$. Hence
	there is useful cancellation and \[ \sum_{\nu \in OP(n)}
	\left( 2 - \frac{8m_1(\mu)}{n} + \frac{6 m_1(\nu)}{n} \right)
	\frac{z_{\nu} p_2(\tilde{K}_{\nu})^2}{n!  2^{n-l(\nu)}} =
	O(n^{-4i-4}). \] It follows that the second term in the upper
	bound of Theorem \ref{steinbound} is $O(n^{-1/4})$, completing
	the proof. \end{proof}

\section{Association Schemes} \label{assoc}

	This final section adapts techniques from earlier sections to
	study the spectrum of an adjacency matrix of an
	association scheme. As is well known (see for instance Chapter
	3 of \cite{BI}) this includes the spectrum a distance regular
	graph as a special case.

	Subsection \ref{rtassoc} discusses needed facts about
	association schemes. Subsection \ref{assocCLT} derives a
	general central limit theorem for the spectrum of an
	association scheme. Subsection \ref{Hamming} illustrates the
	theory of Subsection \ref{assocCLT} for the special case of
	the Hamming scheme, where the result amounts to a central
	limit theorem for the spectrum of the Hamming graph, or
	equivalently for values of q-Krawtchouk polynomials. For an
	application of tools in this section to a problem with a
	non-normal limit, see \cite{CF}.

\subsection{Background on association schemes} \label{rtassoc}

	This subsection gives preliminaries about association
	schemes, using notation from \cite{MS}. Another useful
	reference is \cite{BI}, and what we call association schemes
	some authors call symmetric association schemes.

	{\bf Definition}. An association scheme with $n$ classes
	consists of a finite set $X$ with $n+1$
	relations $R_0,\cdots,R_n$ defined on $X$ which satisfy:

\begin{enumerate}
\item Each $R_i$ is symmetric: $(x,y) \in R_i \Rightarrow (y,x) \in R_i$.
\item For every $x,y \in X$, one has that $(x,y) \in R_i$ for exactly one $i$.
\item $R_0=\{(x,x):x \in X\}$ is the identity relation.
\item If $(x,y) \in R_k$, the number of $z \in X$ such that $(x,z) \in R_i$ and $(y,z) \in R_j$ is a constant $c_{ijk}$ depending on $i,j,k$ but not on the particular choice of $x$ and $y$. 
\end{enumerate}

	The adjacency matrix $D_i$ corresponding to the relation $R_i$
	is the $|X| \times |X|$ matrix with rows and columns labeled
	by the points of $X$ defined by \[ (D_i)_{x,y} = \left\{
	\begin{array}{ll} 1 & \mbox{if $(x,y) \in R_i$}\\ 0 &
	\mbox{otherwise.} \end{array} \right.\] The Bose-Melner
	algebra is defined to be the vector space consisting of all
	matrices $\sum_{i=0}^n a_i D_i$ with $a_i$ real. Since the
	matrices in the Bose-Melner algebra are symmetric and commute
	with each other (by parts 1 and 4 of the definition of an
	association scheme), the Bose-Melner algebra is semisimple and
	has a distinguished basis of primitive idempotents
	$J_0,\cdots,J_n$ satisfying \begin{enumerate} \item
	$J_i^2=J_i$ \ \ \ $i=0,\cdots,n$.  \item $J_i J_k=0$ \ \ \ $i
	\neq k$.  \item $\sum_{i=0}^n J_i = I.$ \end{enumerate} Here I
	is the identity matrix.

 	The $D$'s are also a basis of the Bose-Melner algebra, so one
	can write \[ D_s = \sum_{i=0}^n \phi_s(i) J_i, \ \ \ s = 0,
	\cdots,n\] where the $\phi_s(i)$ are real numbers. Since $D_s
	J_i = \phi_s(i) J_i$, the $\phi_s(i)$ are the eigenvalues of
	$D_s$. For later use note that since $D_0$
	and $\sum_{i=0}^n J_i$ are both the identity matrix,
	$\phi_0(j)=1$ for all $j$.

	 Let $\mu_i$ be the rank of $J_i$. For all $s$, this is the
	multiplicity of the eigenvalue $\phi_s(i)$ of $D_s$. We define
	the Plancherel measure of the association scheme to be the
	probability measure on $\{0,\cdots,n \}$ which chooses $i$
	with probability $\frac{\mu_i}{|X|}$.
	
	Lemmas \ref{asorthog1} and \ref{asorthog2} are the
	orthogonality relations for eigenvalues of association
	schemes. To state them it is helpful to define
	$v_i=c_{ii0}$, so that for any $x$, one has that $v_i$ is the
	number of pairs $(x,y) \in R_i$. Clearly $\sum_{i=0}^n v_i =
	|X|$.

\begin{lemma}
 \label{asorthog1} (\cite{MS}, page 655) For $0 \leq k,l \leq n$, \[
 \sum_{r=0}^n \frac{\phi_r(k) \phi_r(l)}{v_r} = \frac{|X|}{\mu_k}
 \delta_{k,l}.\]
\end{lemma}

\begin{lemma} \label{asorthog2} (\cite{MS}, page 655)  For $0 \leq k,l \leq n$, \[ \sum_{i=0}^n \mu_i \phi_k(i) \phi_l(i) = |X| v_k \delta_{k,l}.\]
\end{lemma}

	Lemma \ref{ascountsol} will be crucial.

\begin{lemma}
 \label{ascountsol} The coefficient of $D_l$ in $D_{s_1} \cdots
 D_{s_m}$ is equal to \[ \frac{1}{v_l} \sum_{i=0}^n \frac{\mu_i}{|X|} \phi_{s_1}(i)
 \cdots \phi_{s_m}(i) \phi_l(i).\] In particular, the coefficient of $D_l$ in
 $(D_s)^m$ is \[ \frac{1}{v_l} \ee[ \phi_s(i)^m \phi_l(i)], \] where $s$ is
 random from Plancherel measure of the association scheme.
\end{lemma}

\begin{proof} Since $D_s = \sum_{i=0}^n \phi_s(i) J_i$, one knows that \[ D_{s_1} \cdots D_{s_m} = \sum_{i=0}^n \phi_{s_1}(i) \cdots \phi_{s_m}(i) J_i.\] By pages 654 and 655 of \cite{MS}, \[ J_i = \frac{1}{|X|} \sum_{l=0}^n \frac{\mu_i \phi_l(i)}{v_l} D_l.\] The result follows. \end{proof}

	Let $p_m(r)$ denote $\frac{v_r}{(v_s)^m}$ multiplied by the
	coefficient of $D_r$ in $(D_s)^m$. It follows from the
	definitions that $p_m(r)$ admits the following probabilistic
	interpretation. Start from some point $x_0 \in X$, move to a
	random $x_1 \in X$ such that $(x_0,x_1) \in R_s$, then to a
	random $x_2 \in X$ such that $(x_1,x_2) \in R_s$, and so on
	until one obtains $x_m$. Then $p_m(r)$ is the probability that
	$(x_0,x_m) \in R_r$.

	The following fact will be useful.

\begin{lemma} \label{asrel}
\begin{enumerate}
\item $p_4(0) = \sum_{r=0}^n \frac{p_2(r)^2}{v_r}$.
\item $p_3(s) = v_s  \sum_{r=0}^n \frac{p_2(r)^2}{v_r}$. 
\end{enumerate}
\end{lemma}

\begin{proof} The first assertion is clear from the probabilistic interpretation of
$p_m(r)$. For an analytic proof, note that by Lemma \ref{ascountsol} (in
the first and fourth equalities) and Lemma \ref{asorthog1} (in the
second and third equalities) \begin{eqnarray*} \sum_{r=0}^n
\frac{p_2(r)^2}{v_r} & = & \sum_{r=0}^n \frac{1}{v_r}
\left( \frac{1}{(v_s)^2} \sum_{i=0}^n \frac{\mu_i}{|X|} \phi_s(i)^2 \phi_r(i) \right)^2\\
& = & \frac{1}{(v_s)^4} \sum_{r=0}^n \frac{1}{v_r} \sum_{i=0}^n \frac{(\mu_i)^2}{|X|^2}
\phi_s(i)^4 \phi_r(i)^2\\ & = & \frac{1}{(v_s)^4} \sum_{i=0}^n \frac{\mu_i}{|X|}
\phi_s(i)^4 \\ & = & p_4(0). \end{eqnarray*}

	The second assertion follows from the first assertion since
$p_3(s)=v_s p_4(0)$, as can be seen either by the probabilistic
interpretation of $p_m(r)$ or by computing both sides of the equation
using Lemma \ref{ascountsol}. \end{proof}

\subsection{Central limit theorem for the spectrum} \label{assocCLT}

	Recall that the goal is to study Plancherel measure of
	association schemes, which chooses $i \in \{0,\cdots,n\}$ with
	probability $\frac{\mu_i}{|X|}$. More precisely, for $s$
	fixed, a central limit theorem is proved for the random
	variable $W$ whose value at $i \in \{0,\cdots,n \}$ is
	$\frac{\phi_s(i)}{\sqrt{v_s}}$.

 	Given $t \in \{0,\cdots,n\}$, we define a Markov chain $L_t$ on
	the set $\{0,\cdots,n \}$ which moves from $i$ to $j$ with
	probability \[ L_t(i,j) := \frac{\mu_j}{|X|} \sum_{r=0}^n
	\frac{\phi_r(i) \phi_r(t) \phi_r(j)}{(v_r)^2} .\]

\begin{lemma} \label{asniceprop} The transition probabilities of $L_t$ are real and non-negative and sum to 1. Moreover the chain $L_t$ is reversible with respect to Plancherel measure.
\end{lemma}

\begin{proof} By Theorems 3.6 and 3.8 of \cite{BI}, the $L_t(i,j)$ are non-negative real numbers. Next, observe that \[ \sum_{j=0}^n L_t(i,j)
 = \frac{1}{|X|} \sum_{r=0}^n \frac{\phi_r(i) \phi_r(t)}{(v_r)^2}
 \sum_{j=0}^n \mu_j \phi_r(j).\] Since $\phi_0(j)=1$ for all $j$,
 Lemma \ref{asorthog2} implies that \[ \sum_{j=0}^n L_t(i,j) = v_0
 \sum_{r=0}^n \frac{\phi_r(i) \phi_r(t)}{(v_r)^2} \delta_{r,0} = 1.\]
 For the reversibility assertion, it is clear from the definition of
 $L_t$ that \[ \frac{\mu_i}{|X|} L_t(i,j) = \frac{\mu_j}{|X|} L_t(j,i)
 \] for all $i,j$. \end{proof}
	
	One uses the chain $L_t$ to construct an exchangeable pair
	$(W,W')$ in the usual way. First choose $i$ from Plancherel
	measure, then choose $j$ with probability $L_t(i,j)$, and
	finally let $(W,W') = (W(i),W(j))$.

\begin{lemma} \label{assteinsat} $\ee(W'|W) = \left( \frac{\phi_s(t)}{v_s} \right) W$.
\end{lemma}

\begin{proof} By the definitions and Lemma \ref{asorthog2}, one has that \begin{eqnarray*} \ee(W'|i) & = & \frac{1}{\sqrt{v_s}} \sum_{j=0}^n
L_t(i,j) \phi_s(j)\\ & = & \frac{1}{\sqrt{v_s}} \frac{1}{|X|}
\sum_{r=0}^n \frac{\phi_r(i) \phi_r(t)}{(v_r)^2} \sum_{j=0}^n \mu_j
\phi_s(j) \phi_r(j)\\ & = & \left( \frac{\phi_s(t)}{v_s} \right)
W(i). \end{eqnarray*} The result follows since this depends on
$i$ only through $W$. \end{proof}

	Corollary \ref{eiglast} will not be used but is worth recording.

\begin{cor} \label{eiglast} The eigenvalues of $L_t$ are $\frac{\phi_s(t)}{v_s}$ for $0 \leq s \leq n$. The functions $\psi_s(i) = \frac{\phi_s(i)}{\sqrt{v_s}}$ are a basis of eigenvectors of $L_t$, orthonormal with respect to the inner product \[ \langle f_1,f_2 \rangle = \sum_{i=0}^n f_1(i) \overline{f_2(i)} \frac{\mu_i}{|X|}.\]
\end{cor}

\begin{proof} The proof of Lemma \ref{assteinsat} shows that $\psi_s$ is an eigenvector of $L_t$ with eigenvalue  $\frac{\phi_s(t)}{v_s}$. The orthonormality assertion follows from Lemma \ref{asorthog2}, and the basis assertion follows since there are $n+1$ eigenvectors. \end{proof}

\begin{lemma} \label{asimmcor} $\ee(W'-W)^2 = 2 \left( 1-\frac{\phi_s(t)}{v_s} \right)$. \end{lemma}

\begin{proof} This is immediate from Lemmas \ref{var} and \ref{assteinsat}.
\end{proof}

\begin{lemma} \label{aspre2} $\ee((W')^2|i) = v_s \sum_{r=0}^n \frac{\phi_r(i)\phi_r(t)}{(v_r)^2} p_2(r)$.
\end{lemma}

\begin{proof} It follows from the definitions that \begin{eqnarray*}
\ee((W')^2|i) & = & \frac{1}{v_s} \sum_{j=0}^n \frac{\mu_j}{|X|}
\sum_{r=0}^n \frac{\phi_r(i) \phi_r(t) \phi_r(j)}{(v_r)^2}
\phi_s(j)^2\\ & = & \frac{1}{v_s} \sum_{r=0}^n \frac{\phi_r(i)
\phi_r(t)}{(v_r)^2} \sum_{j=0}^n \frac{\mu_j}{|X|} \phi_s(j)^2
\phi_r(j). \end{eqnarray*} The result now follows from Lemma
\ref{ascountsol}. \end{proof}

	Note that in the next lemma, the sum is over non-zero $r$. 

\begin{lemma} \label{asbig1} \[ Var(\ee[(W'-W)^2|i]) = (v_s)^2 \sum_{r=1}^n \frac{p_2(r)^2}{v_r} \left( \frac{\phi_r(t)}{v_r} + 1 - \frac{2 \phi_s(t)}{v_s} \right)^2.\]
\end{lemma}

\begin{proof} By Lemmas \ref{assteinsat}, \ref{asimmcor}, and \ref{aspre2}, $Var(\ee[(W'-W)^2|i])$ is equal to 
\begin{eqnarray*}
& & \ee(\ee((W'-W)^2|i)^2) - 4 \left(1 - \frac{\phi_s(t)}{v_s} \right)^2\\
& = & \ee \left[ v_s \sum_{r=0}^n \frac{\phi_r(i) \phi_r(t)}{(v_r)^2} p_2(r) + \left( 1 - \frac{2 \phi_s(t)}{v_s} \right) W^2 \right]^2\\
& &  - 4 \left(1 - \frac{\phi_s(t)}{v_s} \right)^2\\
& = & T_1 + T_2 + T_3 - 4 \left(1 - \frac{\phi_s(t)}{v_s} \right)^2,
\end{eqnarray*} where 
\[ T_1 = (v_s)^2 \ee \left( \sum_{r=0}^n \frac{\phi_r(i) \phi_r(t)}{(v_r)^2} p_2(r) \right)^2 \]
\[ T_2 = 2 \left(1-\frac{2 \phi_s(t)}{v_s} \right) \sum_{r=0}^n
 \frac{\phi_r(t)}{(v_r)^2} p_2(r) \ee \left( \phi_r(i)
 \phi_s(i)^2  \right) \]
\[ T_3 = \left( 1 - \frac{2 \phi_s(t)}{v_s} \right)^2 \ee \left( \frac{\phi_s(i)^4}{(v_s)^2} \right).\]

	Next we simplify these terms. By Lemma \ref{asorthog2}, \[ T_1
	= (v_s)^2 \sum_{r=0}^n \frac{\phi_r(t)^2}{(v_r)^3} p_2(r)^2.\]
	Lemma \ref{ascountsol} implies that \[ T_2 = 2 (v_s)^2 \left(1
	- \frac{2 \phi_s(t)}{v_s} \right) \sum_{r=0}^n
	\frac{\phi_r(t)}{(v_r)^2} p_2(r)^2.\] Lemma \ref{ascountsol}
	and part 1 of Lemma \ref{asrel} imply that \[ T_3
	= (v_s)^2 \left( 1-\frac{2 \phi_s(t)}{v_s} \right)^2 p_4(0) =
	(v_s)^2 \left( 1-\frac{2 \phi_s(t)}{v_s} \right)^2
	\sum_{r=0}^n \frac{p_2(r)^2}{v_r} .\] Thus \[ T_1+T_2+T_3 =
	(v_s)^2 \sum_{r=0}^n \frac{p_2(r)^2}{v_r}
	\left(\frac{\phi_r(t)}{v_r} + 1 - \frac{2 \phi_s(t)}{v_s}
	\right)^2.\] Since $p_2(0)=\frac{1}{v_s}$, the $r=0$ term of
	$T_1+T_2+T_3$ is equal to $4 \left(1 - \frac{\phi_s(t)}{v_s}
	\right)^2$, proving the result. \end{proof}

\begin{lemma} \label{asmom1} Let $k$ be a positive integer.
\begin{enumerate}
\item $\ee(W'-W)^k = (v_s)^{k/2} \sum_{m=0}^k (-1)^{k-m} {k \choose m} \sum_{r=0}^n \frac{\phi_r(t)}{v_r} \frac{p_m(r) p_{k-m}(r)}{v_r}$.
\item $\ee(W'-W)^4 =  v_s^2 \left[ \sum_{r=0}^n \left(8 \left(
1-\frac{\phi_s(t)}{v_s} \right) - 6  \left(1 - \frac{\phi_r(t)}{v_r} \right) \right) \frac{p_2(r)^2}{v_r} \right]$.
\end{enumerate}
\end{lemma}

\begin{proof} For the first assertion, note that
\begin{eqnarray*}
& & \ee((W'-W)^k|i)\\ & = & \frac{1}{(v_s)^{k/2}} \sum_{j=0}^n
\frac{\mu_j}{|X|} \sum_{r=0}^n \frac{\phi_r(i) \phi_r(t)
\phi_r(j)}{(v_r)^2} (\phi_s(j)-\phi_s(i))^k\\ & = &
\frac{1}{(v_s)^{k/2}} \sum_{m=0}^k (-1)^{k-m} {k \choose m}
\phi_s(i)^{k-m} \sum_{r=0}^n \frac{\phi_r(i) \phi_r(t)}{(v_r)^2}\\ & &
\cdot \sum_{j=0}^n \frac{\mu_j}{|X|} \phi_s(j)^m \phi_r(j)\\ & = &
\frac{1}{(v_s)^{k/2}} \sum_{m=0}^k (-1)^{k-m} {k \choose m}
\phi_s(i)^{k-m} (v_s)^m \sum_{r=0}^n \frac{\phi_r(i) \phi_r(t)
p_m(r)}{(v_r)^2}, \end{eqnarray*} where the last equality is Lemma
\ref{ascountsol}. Consequently, \begin{eqnarray*} \ee((W'-W)^k) & = &
\ee (\ee((W'-W)^k|i))\\ & = & (v_s)^{k/2} \sum_{m=0}^k (-1)^{k-m} {k
\choose m}\\ & & \cdot \sum_{r=0}^n \frac{\phi_r(t)}{(v_r)^2} p_m(r)
\sum_{i=0}^n \frac{\mu_i}{|X|} \frac{\phi_s(i)^{k-m}
\phi_r(i)}{(v_s)^{k-m}}. \end{eqnarray*} The first assertion now
follows from Lemma \ref{ascountsol}.

	For the second part, the first assertion gives that \[
\ee(W'-W)^4 = (v_s)^2 \sum_{m=0}^4 (-1)^{m} {4 \choose m} \sum_{r=0}^n
\frac{\phi_r(t)}{v_r} \frac{p_m(r) p_{4-m}(r)}{v_r}.\]
By page 654 of \cite{MS}, $\phi_r(0)=v_r$ for all $r$. Thus
specializing to $t=0$ shows that \[ 0 = (v_s)^2 \sum_{m=0}^4 (-1)^{m}
{4 \choose m} \sum_{r=0}^n \frac{p_m(r) p_{4-m}(r)}{v_r}.\] So for
general $t$, one has that \[ \ee(W'-W)^4 = -(v_s)^2 \sum_{m=0}^4
(-1)^{m} {4 \choose m} \sum_{r=0}^n \left( 1 -\frac{\phi_r(t)}{v_r}
\right) \frac{p_m(r) p_{4-m}(r)}{v_r}.\] The contribution from the
$m=0,4$ terms is 0, since $p_0(r) =0$ if $r \neq 0$. The contribution
from the $m=2$ term is \[ - 6 (v_s)^2 \sum_{r=0}^n \left(1 -
\frac{\phi_r(t)}{v_r} \right) \frac{p_2(r)^2}{v_r}.\] The $m=1,3$
terms are equal and together contribute \begin{eqnarray*} 8 (v_s)^2
\sum_{r=0}^n \left( 1 - \frac{\phi_r(t)}{v_r} \right) \frac{p_1(r)
p_3(r)}{v_r} & = & 8 (v_s)^2 \left( 1 - \frac{\phi_s(t)}{v_s} \right)
\frac{p_3(s)}{v_s}\\ & = & 8 (v_s)^2 \left(1 - \frac{\phi_s(t)}{v_s}
\right) \sum_{r=0}^n \frac{p_2(r)^2}{v_r}, \end{eqnarray*} where the
final equality is part 2 of Lemma \ref{asrel}. Adding the terms
completes the proof of the second assertion. \end{proof}

	Arguing as in the proof of Theorem \ref{main1}, and using the
	above lemmas, one obtains the following result.

\begin{theorem} \label{asmains} Fix $s \in \{0,\cdots,n\}$, and let $W = \frac{\phi_s(i)}{\sqrt{v_s}}$ where $i$ is chosen from the Plancherel measure of the association scheme. Fix $t \in \{1, \cdots, n\}$. Then for all real $x_0$,
\begin{eqnarray*}
& & \left| \pp(W \leq x_0) - \frac{1}{\sqrt{2 \pi}}
\int_{-\infty}^{x_0} e^{-\frac{x^2}{2}} dx \right|\\ & \leq &
\frac{v_s}{a} \sqrt{ \sum_{r=1}^n \frac{p_2(r)^2}{v_r} \left(
\frac{\phi_r(t)}{v_r}+1 - \frac{2 \phi_s(t)}{v_s} \right)^2}\\
 & & +
\frac{\sqrt{v_s}}{(\pi)^{1/4}} \left[ \sum_{r=0}^n \left( 8
-\frac{6}{a} \left(1-\frac{\phi_r(t)}{v_r} \right) \right)
\frac{p_2(r)^2}{v_r} \right]^{1/4}, \end{eqnarray*} where $a=1
-\frac{\phi_s(t)}{v_s}$.
\end{theorem}

\subsection{Example: Hamming Scheme} \label{Hamming}

	This subsection illustrates the theory of Subsection
	\ref{assocCLT} for the Hamming scheme $H(d,q)$, where $d,q$
	are positive integers and $q \geq 2$.

	To begin we recall the definition of $H(d,q)$ and its basic
	properties, referring the reader to Chapter 3 of \cite{BI} for
	more details. The elements $X$ of $H(d,q)$ are d-tuples of
	numbers chosen from $\{1,\cdots,q\}$; clearly $|X|=q^d$. A
	pair $(x,y)$ is in $R_i$ if $x$ and $y$ differ in exactly i
	coordinates. For $0 \leq i \leq d$, one has that $v_i=(q-1)^i
	{d \choose i}$ and $\mu_i=(q-1)^i {d \choose i}$. Thus the
	Plancherel measure on $\{0,\cdots, d\}$ chooses $i$ with
	probability $\frac{(q-1)^i {d \choose i}}{q^d}$. The numbers
	$\phi_s(i)$ are equal to \[ \sum_{j=0}^s (-1)^j (q-1)^{s-j} {i
	\choose j} {d-i \choose s-j}.\] The polynomial $\phi_s$ in the
	variable $i$ is known as a q-Krawtchouk polynomial.

	In what follows $W=\frac{\phi_1(i)}{\sqrt{v_1}}$, where $i$ is
 chosen from Plancherel measure of the Hamming scheme $H(d,q)$. The
 Hamming graph has vertex set $X$ and an edge between two vertices if
 they differ in one coordinate. Thus the eigenvalues of the adjacency
 matrix of the Hamming graph are $\phi_1(i)$ with multiplicity
 $\mu_i$, which motivates the study of $W$. Hora \cite{Ho2} shows
 that if $d \rightarrow \infty$ and $q/d \rightarrow 0$, then $W$
 converges in distribution to a normal random variable with mean 0 and
 variance 1.
	
	In fact since $W(i)=\frac{(q-1)d-qi}{\sqrt{(q-1)d}}$ and
	$\frac{\mu_i}{|X|} = \frac{(q-1)^i {d \choose i}}{q^d}$, it is straightforward
	that $W$ has the same distribution as \[
	\frac{Y_1+\cdots+Y_d}{\sqrt{(q-1)d}}, \] where the $Y$'s are
	independent random variables, each equal to $q-1$ with
	probability $\frac{1}{q}$ and to $-1$ with probability $1-\frac{1}{q}$. Hence
	the Berry-Esseen theorem \cite{Du} shows that $W$ satisfies a
	central limit theorem with the error term $C
	\sqrt{\frac{q}{d}}$ where $C$ is a small explicit constant.

	For comparison, let us study $W$ using Theorem \ref{asmains}
	with $s=1$ and $t=1$. Then $a=\frac{q}{(q-1)d}$ and one
	computes that $p_2(0)=\frac{1}{(q-1)d}$, $p_2(1) =
	\frac{(q-2)}{(q-1)d}$, and $p_2(2)=1-\frac{1}{d}$. The first
	term in the error term of Theorem \ref{asmains} is then
	computed to equal $\sqrt{\frac{(q-2)^2}{(q-1)d}}$, which is
	less than $\sqrt{\frac{q}{d}}$. The second term in the error
	term of Theorem \ref{asmains} is computed to equal $\left(
	\frac{2q^2}{\pi (q-1)d} \right)^{1/4}$, which is at most
	$\left( \frac{2q}{d} \right)^{1/4}$. Thus one obtains a
	central limit theorem for $W$ with error term
	$\sqrt{\frac{q}{d}} +\left( \frac{2q}{d} \right)^{1/4}$.

	To close this section, it is shown how our exchangeable pair
	can be used to give $O(d^{-1/2})$ bounds for $q$ fixed. The
	key is Lemma \ref{bdeath}, which shows that the Markov chain
	$L_1$ is actually a birth-chain and computes its transition
	probabilities.

\begin{lemma} \label{bdeath} The Markov chain $L_1$ on the set $\{0,\cdots,d\}$ is a birth-death chain with transition probabilities \[ L_1(i,j) = \left\{ \begin{array}{ll} 
\frac{i}{d(q-1)} & \mbox{if $j=i-1$}\\
\frac{i}{d} \left(1-\frac{1}{q-1} \right) & \mbox{if $j=i$}\\
1-\frac{i}{d} & \mbox{if $j=i+1$} \end{array} \right. \]
\end{lemma}	

\begin{proof} Recall that $L_1(i,j) = \frac{\mu_j}{|X|} \sum_{r=0}^n \frac{\phi_r(i) \phi_r(1) \phi_r(j)}{(v_r)^2}$. From the formula for $\phi_r(l)$ one sees that \[ \phi_r(l) =  (q-1)^{r-l} \frac{{d \choose r}}{{d \choose l}} \phi_l(r)  \] for any $0 \leq r,l \leq d$. Thus \[ L_1(i,j) = \frac{\mu_j}{|X| d {d \choose i} {d \choose j} (q-1)^{i+j+1}} \sum_{r=0}^n {d \choose r} (q-1)^r \phi_1(r) \phi_i(r) \phi_j(r).\] From page 152 of \cite{MS}, there is a recurrence relation \[
	(i+1) \phi_{i+1}(r) = [(d-i)(q-1)+i-qr] \phi_i(r) -
	(q-1)(d-i+1) \phi_{i-1}(r) .\] Since $\phi_1(r)= d(q-1)-qr$,
	the recurrence is equivalent to \[ \phi_1(r) \phi_i(r) =
	i(q-2) \phi_i(r) + (i+1) \phi_{i+1}(r) + (q-1)(d-i+1)
	\phi_{i-1}(r).\] Applying this to expression for $L_1(i,j)$ at
	the end of the previous paragraph, the result follows from
	Lemma \ref{asorthog2}. \end{proof}

	Lemma \ref{bdeath} implies that $|W'-W| \leq A$ with
$A=\frac{q}{\sqrt{(q-1)d}}$. Thus one can apply the version of Theorem
\ref{asmains} which would arise by using Theorem \ref{rinrot} instead
of Theorem \ref{steinbound}. Recalling that $a=\frac{q}{(q-1)d}$, one
obtains an error term of $\sqrt{\frac{q}{d}} + \frac{.41q^2 + 1.5
q}{\sqrt{(q-1)d}}$. For $q$ fixed this goes as $d^{-1/2}$.

\section{Acknowledgements} The author was supported by NSA grant H98230-05-1-0031 and NSF grant DMS-0503901.

\end{document}